\documentclass{hjm1}
\usepackage{verbatim}

\usepackage{setspace}

\begin{document}

\title[Effective $H^{\infty}$ interpolation]
{Effective $H^{\infty}$ interpolation}
\author[Zarouf Rachid] {Rachid Zarouf}
\address [] {CMI-LATP, UMR 6632, Universit\'e de Provence, 39, rue F.-Joliot-Curie,
13453 Marseille cedex 13, France} \email{rzarouf@cmi.univ-mrs.fr}\

\thanks{I would like to thank Professor Nikolai Nikolski for all of his work,
his wisdom and the pleasure that our discussions gave to me.  I also would like to thank the referee for the careful
review and the valuable comments, which provided insights that helped
improve the paper.}



\begin{abstract}
Given a finite subset $\sigma$ of the unit disc $\mathbb{D}$ and
a holomorphic function $f$ in $\mathbb{D}$ belonging to a class
$X$, we are looking for a function $g$ in another class $Y$ which
satisfies $g_{|\sigma}=f_{|\sigma}$ and is of minimal norm
in $Y$. Then, we wish to compare $\left\Vert g\right\Vert _{Y}$
with $\left\Vert f\right\Vert _{X}.$ More precisely, we consider
the interpolation constant $c\left(\sigma,\, X,\, Y\right)=\mbox{sup}{}_{f\in X,\,\parallel f\parallel_{X}\leq1}\mbox{inf}_{g_{|\sigma}=f_{|\sigma}}\left\Vert g\right\Vert _{Y}.$
When $Y=H^{\infty}$, our interpolation problem includes those of
Nevanlinna-Pick and Caratheodory-Schur.  Moreover, Carleson's free
interpolation problem can be interpreted in terms of the constant
$c\left(\sigma,\, X,\, H^{\infty}\right)$. For $Y=H^{\infty}$, $X=H^{p}$
(the Hardy space) or $X=L_{a}^{2}$ (the Bergman space), we obtain
an upper bound for the constant $c\left(\sigma,\, X,\, H^{\infty}\right)$
in terms of $n=\mbox{card}\,\sigma$ and $r=\mbox{max}{}_{\lambda\in\sigma}\left|\lambda\right|$.
Our upper estimates are shown to be sharp with respect to $n$ and
$r$.  \end{abstract}

\maketitle


\newtheorem*{Prov}{Provision}
\newtheorem {mythm}{Theorem}[section]
\newtheorem {mylem}{Lemma}
\newtheorem*{idfn}{Definition}
\newenvironment{definition}{\begin{idfn}
\rm}{\end{idfn}}

\section{Introduction}\label{section1}

\subsection{Statement and historical context of the problem}\label{subsection11}

Let $\mathbb{D}=\{z\in\mathbb{C}:\,\vert z\vert<1\}$ be the unit
disc of the complex plane and let ${\rm Hol}\left(\mathbb{D}\right)$
be the space of holomorphic functions on $\mathbb{D}.$ The problem
considered is the following: given two Banach spaces $X$ and $Y$
of holomorphic functions on $\mathbb{D},$ $X,\, Y\subset{\rm Hol}\left(\mathbb{D}\right),$
and a finite subset $\sigma\subset\mathbb{D}$, find the least norm
interpolation by functions of the space $Y$ for the traces $f_{\vert\sigma}$
of functions of the space $X$, in the worst case of $f$. The case
$X\subset Y$ is of no interest, and so one can suppose that either
$Y\subset X$ or $X,\, Y$ are incomparable.

More precisely, our problem is to compute or estimate the following
interpolation constant
\[
c\left(\sigma,\, X,\, Y\right)=\sup_{f\in X,\,\parallel f\parallel_{X}\leq1}\inf\left\{ \left\Vert g\right\Vert _{Y}:\, g_{|\sigma}=f_{|\sigma}\right\} .\]
If $r\in[0,\,1)$ and $n\geq1,$ we also define
\[
C_{n,\, r}(X,\, Y)=\sup\left\{ c(\sigma,\, X,\, Y)\,:\;{\rm card\,}\sigma\leq n\,,\,\left|\lambda\right|\leq r,\;\forall\lambda\in\sigma\right\} .\]
Here and later on, $H^{\infty}$ stands for the space (algebra) of
bounded holomorphic functions on $\mathbb{D}$ endowed with the norm
$\left\Vert f\right\Vert _{\infty}=\sup_{z\in\mathbb{D}}\left|f(z)\right|.$
The classical interpolation problems -those of Nevanlinna-Pick{\small{}
}(1916) and Carathéodory-Schur (1908) (see \cite{Nik2}
 p.231 for these
two problems) on the one hand, and Carleson's free interpolation (1958)
(see {[}16{]} p.158) on the other hand- are of this nature and correspond
to the case $Y=H^{\infty}$. Two first are {}``individual'', in
the sense that one looks simply to compute the norms $\left\Vert f\right\Vert _{H_{|\sigma}^{\infty}}$
or $\left\Vert f\right\Vert _{H^{\infty}/z^{n}H^{\infty}}$ for a
given $f$. In the case of the third one, we consider infinite sets
$\sigma.$ Let $l^{\infty}(\sigma)$ be the space of bounded functions
$\left(a_{\lambda}\right)_{\lambda\in\sigma}$ on $\sigma$ endowed
with the norm $\left\Vert a\right\Vert _{l^{\infty}(\sigma)}=\max_{\lambda\in\sigma}\left|a_{\lambda}\right|$.
Carleson's free interpolation problem is to compare the norms $\left\Vert a\right\Vert _{l^{\infty}(\sigma)}$
and \[
\mbox{inf}\left\{ \left\Vert g\right\Vert _{\infty}:\, g(\lambda)=a_{\lambda},\:\lambda\in\sigma\right\} .\]
In other words, we want to estimate the interpolation constant defined
as
\[
c\left(\sigma,\, l^{\infty}(\sigma),\, H^{\infty}\right)=\sup_{a\in l^{\infty}(\sigma),\,\left\Vert a\right\Vert _{l^{\infty}}\leq1}\mbox{inf}\left\{ \left\Vert g\right\Vert _{\infty}:\, g(\lambda)=a_{\lambda},\:\lambda\in\sigma\right\} .\]
Let us now explain why our problem includes those of Nevanlinna-Pick
and Carathéodory-Schur.

\textbf{$\ \ \ \ \ $}(i) Nevannlinna-Pick interpolation problem.

Given $\sigma=\left\{ \lambda_{1},\,...,\,\lambda_{n}\right\} $ a
finite subset of $\mathbb{D}$ and $\mathcal{W}=\left\{ w_{1},\,...,\, w_{n}\right\} $
a finite subset of $\mathbb{C},$ find\[
{\rm NP}{}_{\sigma,\,\mathcal{W}}=\inf\left\{ \left\Vert f\right\Vert _{\infty}:\, f\left(\lambda_{i}\right)=w_{i},\, i=1,\,...,\, n\right\} .\]
The classical answer of Pick is the following:\[
{\rm NP}{}_{\sigma,\,\mathcal{W}}=\inf\left\{ c>0:\,\left(\frac{c^{2}-\overline{w_{i}}w_{j}}{1-\overline{\lambda_{i}}\lambda_{j}}\right)_{1\leq i,\, j\leq n}\gg0\right\} ,\]
where for any $n\times n$ matrix $M$, $M\gg0$ means that $M$ is
positive definite.

\textbf{$\ \ \ \ \ $}(ii) Carathéodory-Schur interpolation problem.

Given $\mathcal{W}=\left\{ w_{0},\, w_{1},\,...,\, w_{n}\right\} $
a finite subset of $\mathbb{C},$ find
\[
{\rm CS}_{\mathcal{W}}=\inf\left\{ \left\Vert f\right\Vert _{\infty}:\, f\left(z\right)=w_{0}+w_{1}z+...+w_{n}z^{n}+...\right\} .\]
The classical answer of Schur is the following:\[
{\rm CS}_{\mathcal{W}}=\left\Vert \left(T_{\varphi}\right)_{n}\right\Vert ,\]
where $T_{\varphi}$ is the Toeplitz operator associated with a symbol
$\varphi\,,$ $\left(T_{\varphi}\right)_{n}$ is the compression of
$T_{\varphi}$ on $\mathcal{P}_{n},$ the space of analytic polynomials
of degree less or equal than $n$, and $\varphi$ is the polynomial
$\sum_{k=0}^{n}w_{k}z^{k}$.
From a modern point of view, these two interpolation problems (i)
and (ii) are included in the following mixed problem: given $\sigma=\left\{ \lambda_{1},\,...,\,\lambda_{n}\right\} \subset\mathbb{D}$
and $f\in{\rm Hol}(\mathbb{D})$, compute or estimate
\[
\left\Vert f\right\Vert _{H^{\infty}/B_{\sigma}H^{\infty}}=\inf\left\{ \left\Vert g\right\Vert _{\infty}:\, f-g\in B_{\sigma}{\rm Hol}(\mathbb{D})\right\} .\]
From now on, if $\sigma=\left\{ \lambda_{1},\,...,\,\lambda_{n}\right\} \subset\mathbb{D}$
is a finite subset of the unit disc, then \[
B_{\sigma}=\prod_{j=1}^{n}b_{\lambda_{j}}\]
is the corresponding finite Blaschke product where $b_{\lambda}=\frac{\lambda-z}{1-\overline{\lambda}z},$
$\lambda\in\mathbb{D}$. The classical Nevanlinna-Pick problem corresponds
to the case $X={\rm Hol}(\mathbb{D})$, $Y=H^{\infty}$, and the one
of Carathéodory-Schur to the case $\lambda_{1}=\lambda_{2}=...=\lambda_{n}=0$
and $X={\rm Hol}(\mathbb{D})$, $Y=H^{\infty}$.

Looking at this problem in the form of computing or estimating the
interpolation constant $c\left(\sigma,\, X,\, Y\right)$ which is
nothing but the norm of the embedding operator $\left(X_{|\sigma},\,\left\Vert \cdot\right\Vert _{X_{|\sigma}}\right)\rightarrow\left(Y_{|\sigma},\,\left\Vert \cdot\right\Vert _{Y_{|\sigma}}\right)$,
one can think, of course, on passing (after) to the limit -in the
case of an infinite sequence $\left\{ \lambda_{j}\right\} $ and its
finite sections $\left\{ \lambda_{j}\right\} _{j=1}^{n}$- in order
to obtain a Carleson type interpolation theorem $X_{|\sigma}=Y_{|\sigma},$
but not necessarily. In particular, even the classical Nevanlinna-Pick
theorem (giving a necessary and sufficient condition on a function
$a$ for the existence of $f\in H^{\infty}$ such that $\left\Vert f\right\Vert _{\infty}\leq1$
and $f(\lambda)=a_{\lambda},$ $\lambda\in\sigma$), does not lead
immediately to Carleson's criterion for $H_{|\sigma}^{\infty}=l^{\infty}(\sigma).$
(Finally, a direct deduction of Carleson's theorem from Pick's result
was done by P. Koosis \cite{Koo} in 1999 only). Similarly, the problem
stated for $c\left(\sigma,\, X,\, Y\right)$ is of interest in its
own. It is a kind of {}``effective interpolation'' because we are
looking for sharp estimates or a computation of $c\left(\sigma,\, X,\, Y\right)$
for a variety of norms $\left\Vert \cdot\right\Vert _{X},\;\left\Vert \cdot\right\Vert _{Y}.$

\subsection{Motivations } \label{subsection12}

$\;$

\textbf{a.} As it is mentioned in Subsection \ref{subsection11}, one of the most
interesting cases is $Y=H^{\infty}$. In this case, the quantity $c\left(\sigma,\, X,\, H^{\infty}\right)$
has a meaning of an intermediate interpolation between the Carleson
one (when $\left\Vert f\right\Vert _{X_{\vert\sigma}}\asymp{\displaystyle \sup_{1\leq i\leq n}}\left|f\left(\lambda_{i}\right)\right|$)
and the individual Nevanlinna-Pick interpolation (no conditions on
$f$).

\textbf{b.} The following partial case was especially stimulating
(which is a part of a more complicated question arising in an applied
situation in \cite{Barat, BaratWielon}): given a set $\sigma\subset\mathbb{D}$,
how can one estimate $c\left(\sigma,\, H^{2},\, H^{\infty}\right)$
in terms of $n={\rm card}\,\sigma$ and ${\displaystyle \max_{\lambda\in\sigma}}\left|\lambda\right|=r$
only? (Here, $H^{2}$ is the standard Hardy space of the disc $\mathbb{D}$
and is defined below in Subsection \ref{subsection13}).

\textbf{c.} There is a direct link between the constant $c\left(\sigma,\, X,\, Y\right)$
and numerical analysis. For example, in matrix analysis, it is of
interest to bound the norm of an $H^{\infty}$-calculus $\left\Vert f(A)\right\Vert \leq c\left\Vert f\right\Vert _{\infty}$,
$f\in H^{\infty},$ for a contraction $A$ on an $n$-dimensional
arbitrary Banach space, with a given spectrum $\sigma(A)\subset\sigma$.
The best possible constant is $c=c\left(\sigma,\, H^{\infty},\, W\right)$,
so that
\[
c\left(
\sigma,\, H^{\infty},\, W
\right)=
\]
\[=\sup_{\parallel f\parallel_{\infty}\leq1}\sup\left\{ \left\Vert f(A)\right\Vert :\, A:\left(\mathbb{C}^{n},\,\vert\cdot\vert\right)\rightarrow\left(\mathbb{C}^{n},\,\vert\cdot\vert\right),\,\left\Vert A\right\Vert \leq1,\,\sigma(A)\subset\sigma\right\} ,
\]
where $W=\left\{ f=\sum_{k\geq0}\hat{f}(k)z^{k}:\:\sum_{k\geq0}\left|\hat{f}(k)\right|<\infty\right\} $
stands for the Wiener algebra, and the interior sup is taken over
all contractions on $n-$dimensional Banach spaces.\textbf{ }An estimate
for $c\left(\sigma,\, H^{\infty},\, W\right)$ is given in \cite{Nik3}.
An interesting case occurs for $f$ such that $f_{\vert\sigma}=\frac{1}{z}\vert_{\sigma}$
(estimates on condition numbers and the norm of inverses of $n\times n$
matrices) or $f_{|\sigma}=\frac{1}{\lambda-z}_{|\sigma}$ (for estimates
on the norm of the resolvent of an $n\times n$ matrix). Notice that
in the same spirit, the case $Y=B_{\infty,1}^{0}$ where $B_{\infty,1}^{0}$
is a Besov algebra presents an interesting case for the functional
calculus of finite rank operators, in particular, those satisfying
the so-called Ritt condition.

\subsection{The spaces $X$ and $Y$ considered here}\label{subsection13}

We systematically use the following conditions for the spaces $X$
and $Y$,

\def\theequation{$P_{1}$}
\begin{equation}
{\rm Hol}((1+\epsilon)\mathbb{D})\:{\rm is}\:{\rm continuously\: embedded\: into\:}Y\:{\rm for\: every}\:\epsilon>0,\label{eq:}\end{equation}

\def\theequation{$P_{2}$}
\begin{equation}
Pol_{+}\subset X\:{\rm and}\; Pol_{+}{\rm \: is\: dense\: in}\: X,\label{eq:}\end{equation}
where $Pol_{+}$ stands for the set of all complex polynomials $p$,
$p=\sum_{k=0}^{N}a_{k}z^{k},$
\def\theequation{$P_{3}$}
\begin{equation}
\left[f\in X\right]\Rightarrow\left[z^{n}f\in X\,,\,\forall n\geq0\:{\rm and\:\overline{lim}}\left\Vert z^{n}f\right\Vert ^{\frac{1}{n}}\leq1\right],\label{eq:}\end{equation}
\def\theequation{$P_{4}$}\begin{equation}
\left[f\in X,\,\lambda\in\mathbb{D}\;{\rm and}\; f(\lambda)=0\right]\Rightarrow\left[\frac{f}{z-\lambda}\in X\right].\label{eq:}\end{equation}

Assuming $X$ satisfies property $(P_{4})$ and $Y\subset X$, then
the quantity $c\left(\sigma,\, X,\, Y\right)$ can be written as follows
\[
c\left(\sigma,\, X,\, Y\right)=\sup_{\parallel f\parallel_{X}\leq1}\inf\left\{ \left\Vert g\right\Vert _{Y}:\, g\in Y,\, g-f\in B_{\sigma}X\right\} .\]
General spaces $X$ and $Y$ satisfying $(P_{i})_{1\leq i\leq4}$
are studied in Section \ref{section3}. Then, we study special cases of such spaces:
from Section \ref{section4} to the end of this paper, $Y=H^{\infty},$ but $X$
may change from one section to another. In particular, in Sections
\ref{section4} and \ref{section5}, $X=H^{p}=H^{p}(\mathbb{D}),$ $1\leq p\leq\infty$ which
are the standard Hardy spaces on the disc $\mathbb{D}$ (see {[}15{]}
Chapter 2) of all $f\in{\rm Hol}(\mathbb{D})$ satisfying \[
\sup_{0\leq r<1}\left(\int_{\mathbb{T}}\left|f(rz)\right|^{p}{\rm d}m(z)\right)^{1/p}<\infty,\]
$m$ being the Lebesgue normalized measure on $\mathbb{T}.$ From
now on, if $f\in{\rm Hol}(\mathbb{D})$ and $k\in\mathbb{N},$ \[
\hat{f}(k)\;{\rm stands}\;{\rm for}\;{\rm the}\; k^{th}\;{\rm Taylor}\;{\rm coefficient}\;{\rm of}\; f.\]
 For $p=2,$ an equivalent description of $H^{2}$ is \[
H^{2}=\left\{ f=\sum_{k\geq0}\hat{f}(k)z^{k}:\,\sum_{k\geq0}\left|\hat{f}(k)\right|^{2}<\infty\right\} .\]
We also study (see Section \ref{section6}) the case $X=l_{a}^{2}\left(1/\sqrt{k+1}\right),$
which is the Bergman space of all $f=\sum_{k\geq0}\hat{f}(k)z^{k}$
satisfying \[
\sum_{k\geq0}\left|\hat{f}(k)\right|^{2}\frac{1}{k+1}<\infty.\]
This space is also given by: $X=L_{a}^{2}$, the space of holomorphic
functions $f$ on $\mathbb{D}$ such that \[
\int_{\mathbb{D}}\left|f(z)\right|^{2}{\rm d}A<\infty,\]
where ${\rm d}A$ stands for the area measure.

\section{Results}\label{section2}

We start studying general Banach spaces $X$ and $Y$ and give some
sufficient conditions under which $C_{n,\, r}(X,\, Y)<\infty$. In
particular, we prove the following fact.

\begin{thm} \label{Theorem A}Let $X,\, Y$ be Banach
spaces satisfying properties $\left(P_{i}\right)$, $i=1,\,...,\,4$.
Then
\[
C_{n,\, r}(X,\, Y)<\infty,\]
for every $n\geq1$ and $r\in[0,\,1).$
\end{thm}

Next, we add the condition that $X$ is a Hilbert space, and give
in this case a general upper bound for the quantity $C_{n,\, r}(X,\, Y)$.

\begin{thm} \label{Theorem B}Let $Y$ be a Banach space
satisfying property $\left(P_{1}\right)$ and $X=\left(H,\,\left(\cdot,\,\cdot\right)_{H}\right)$
a Hilbert space satisfying properties $\left(P_{i}\right)$ for $i=2,\,3,\,4$.
We moreover suppose that for every $0<r<1$ there exists $\epsilon>0$
such that $k_{\lambda}\in{\rm Hol}\left((1+\epsilon)\mathbb{D}\right)$
for all $\vert\lambda\vert<r$, where $k_{\lambda}$ stands for the
reproducing kernel of $X$ at point $\lambda$, and $\overline{\lambda}\mapsto k_{\lambda}$
is holomorphic on $\vert\lambda\vert<r$ as a ${\rm Hol}((1+\epsilon)\mathbb{D})$-valued
function. Let $\sigma=\{\lambda_{1},...,\lambda_{1},\lambda_{2},...,\lambda_{2},...,\lambda_{t},...,\lambda_{t}\}$
be a sequence in $\mathbb{D}$, where $\lambda_{s}$ are repeated
according to their multiplicity $m_{s}$, $\sum_{s=1}^{t}m_{s}=n$.
Then we have,

i)\[
c\left(\sigma,\, H,\, Y\right)\leq\left(\sum_{k=1}^{n}\left\Vert \mathcal{E}_{k}\right\Vert _{Y}^{2}\right)^{\frac{1}{2}},\]
 where $\left(\mathcal{E}_{k}\right)_{1\leq k\leq n}$ stands for
the Gram-Schmidt orthogonalization (in the space $H$) of the sequence
\[
k_{\lambda_{1},0},\, k_{\lambda_{1},1},\,...,\,
k_{\lambda_{1},m_{1}-1},\, k_{\lambda_{2},0},\, k_{\lambda_{2},1}\,,...,\,
k_{\lambda_{2},m_{2}-1},...,
\, k_{\lambda_{t},0},\, k_{\lambda_{t},1},\,...,\,
k_{\lambda_{t},m_{t}-1},\]
and $k_{\lambda,\, i}=\left(\frac{d}{d\overline{\lambda}}\right)^{i}k_{\lambda}$,
$i\in\mathbb{N}.$

\vspace{0.3cm}

ii) For the case $Y=H^{\infty}$, we have
\[
c(\sigma,\, H,\, H^{\infty})\leq\sup_{\zeta\in\mathbb{D}}\left\Vert P_{B_{\sigma}}^{H}k_{\zeta}\right\Vert _{H},\]
where $P_{B_{\sigma}}^{H}=\sum_{k=1}^{n}\left(\cdot,\,\mathcal{E}_{k}\right)_{H}\mathcal{E}_{k}$
stands for the orthogonal projection of $H$ onto $K_{B_{\sigma}}(H)$,\[
K_{B_{\sigma}}(H)={\rm span}\left(k_{\lambda_{j},\, i}:\,1\leq i<m_{j},\, j=1,...,\, t\right).\]
\end{thm}

After that, we deal with $H^{\infty}$ interpolation ($Y=H^{\infty})$.
For general Banach spaces (of analytic functions in $\mathbb{D}$)
of moderate growth $X$, we formulate the following conjecture: \[
c_{1}\varphi_{X}\left(1-\frac{1-r}{n}\right)\leq C_{n,r}\left(X,\, H^{\infty}\right)\leq c_{2}\varphi_{X}\left(1-\frac{1-r}{n}\right),\]
where $\varphi_{X}(t)$, $0\leq t<1$ stands for the norm of the evaluation
functional $f\mapsto f(t)$ on the space $X$. We prove this conjecture
for $X=H^{p},\, L_{a}^{2}$, $p\in[1,\,\infty)$ (defined in Subsection \ref{subsection13}). More precisely, we prove the following Theorems \ref{Theorem C} and \ref{Theorem D}. Here
and later on, \[
\sigma_{n,\,\lambda}=\underbrace{\{\lambda,\lambda,...,\lambda\}}_{n},\]
is the one-point set of multiplicity $n$ corresponding to $\lambda\in\mathbb{D}.$

\begin{thm}  \label{Theorem C}  Let $1\leq p\leq\infty$,
$n\geq1$, $r\in[0,\,1),$ and $\left|\lambda\right|\leq r$. We have,\[
\frac{1}{32^{\frac{1}{p}}}\left(\frac{n}{1-\left|\lambda\right|}\right)^{\frac{1}{p}}\leq c\left(\sigma_{n,\,\lambda},\, H^{p},H^{\infty}\right)\leq C_{n,r}\left(H^{p},\, H^{\infty}\right)\leq A_{p}\left(\frac{n}{1-r}\right)^{\frac{1}{p}},\]
where $A_{p}$ is a constant depending only on $p$ and the left-hand
side inequality is proved only for $p\in2\mathbb{Z}_{+}.$ For $p=2$,
we have $A_{2}=\sqrt{2}$. \end{thm}

\begin{thm} \label{Theorem D} Let $1\leq p\leq\infty$,
$n\geq1$, $r\in[0,\,1),$ and $\left|\lambda\right|\leq r$. We have,\[
\frac{1}{32^{\frac{1}{p}}}\left(\frac{n}{1-\left|\lambda\right|}\right)^{\frac{1}{p}}\leq c\left(\sigma_{n,\,\lambda},\, H^{p},H^{\infty}\right)\leq C_{n,r}\left(H^{p},\, H^{\infty}\right)\leq A_{p}\left(\frac{n}{1-r}\right)^{\frac{1}{p}},\]
where $A_{p}$ is a constant depending only on $p$ and the left-hand
side inequality is proved only for $p\in2\mathbb{Z}_{+}.$ For $p=2$,
we have $A_{2}=\sqrt{2}$. \end{thm}

The above Theorems  \ref{Theorem A} , \ref{Theorem C} and \ref{Theorem D}
 were already announced in the note \cite{Z2}.

\vspace{0.1cm}

In order to prove (Theorem \ref{Theorem A} , Theorem \ref{Theorem B}  and) the right-hand side inequality
of Theorem \ref{Theorem C} and Theorem \ref{Theorem D}, given $f\in X$ and $\sigma$ a finite subset
of $\mathbb{D}$, we first use a linear interpolation: \[
f\mapsto\sum_{k=1}^{n}\left\langle f,\, e_{k}\right\rangle e_{k},\]
where $\left\langle \cdot,\,\cdot\right\rangle $ means the Cauchy
sesquilinear form $\left\langle h,\, g\right\rangle =\sum_{k\geq0}\hat{h}(k)\overline{\hat{g}(k)},$
and $\left(e_{k}\right)_{1\leq k\leq n}$ is the Malmquist basis (effectively
constructible) of the space $K_{B}=H^{2}\ominus BH^{2}$, whith $B=B_{\sigma}$
(see N. Nikolski, \cite{Nik3} p. 117)). Next, we use the complex interpolation
between Banach spaces, (see H. Triebel \cite{Tr} Theorem 1.9.3-(a)
p.59). Among the technical tools used in order to find an upper bound
for $\left\Vert \sum_{k=1}^{n}\left\langle f,\, e_{k}\right\rangle e_{k}\right\Vert _{\infty}$
(in terms of $\left\Vert f\right\Vert _{X}$), the most important
one is a Bernstein-type inequality $\left\Vert f'\right\Vert _{p}\leq c_{p}\left\Vert B'\right\Vert _{\infty}\left\Vert f\right\Vert _{p}$
for a (rational) function $f$ in the star-invariant subspace $K_{B}^{p}:=H^{p}\cap B\overline{zH^{p}}$
, $1\leq p\leq\infty$ (for $p=2,$ $K_{B}^{2}=K_{B})$, generated
by a (finite) Blaschke product $B$, (K. Dyakonov \cite{Dya1, Dya2}). For
$p=2$, we give an alternative proof of the Bernstein-type estimate
we need and the constant $c_{2}$ we obtain is slightly better, see
Section \ref{section6}.

The lower bound problem of Theorems \ref{Theorem C} and \ref{Theorem D} is treated in Section
7 by using the {}``worst'' interpolation $n$-tuple $\sigma=\sigma_{n,\,\lambda},$
(the Carathéodory-Schur type interpolation). The {}``worst'' interpolation
data comes from the Dirichlet kernels $\sum_{k=0}^{n-1}z^{k}$ transplanted
from the origin to $\lambda.$ We notice that spaces $X=H^{p},\, L_{a}^{2}$,
$p\in[1,\,\infty)$ satisfy the condition $X\circ b_{\lambda}\subset X$
which makes the problem of upper and lower bound easier.

The paper is organized as follows. In Section \ref{section3} we prove Theorems
\ref{Theorem A} and \ref{Theorem B}. Sections \ref{section4} and \ref{section5} (resp. Section \ref{section6}) are (resp. is) devoted to the proof of the upper estimate of Theorem \ref{Theorem C} (resp. Theorem \ref{Theorem D}).
In Section \ref{section7}, we prove the lower bounds stated in Theorems  \ref{Theorem C} and  \ref{Theorem D}.
At the end of the paper, we shortly compare the method used in Sections
\ref{section3}, \ref{section4}, \ref{section5} and \ref{section6} with those resulting from the Carleson free interpolation,
see Section \ref{section8}.

\section{Upper bounds for $c(\sigma,\, X,\, Y)$}\label{section3}

\subsection{Banach spaces $X$, $Y$ satisfying properties $\left(P_{i}\right)_{1\leq i\leq4}$ }\label{subsection31}

In this subsection, $X$ and $Y$ are Banach spaces which satisfy
properties $\left(P_{i}\right)_{1\leq i\leq4}.$ We prove Theorem \ref{Theorem A}  which shows that in this case our interpolation constant $c(\sigma,X,Y)$
is bounded by a quantity which depends only on $n={\rm card\,}\sigma$
and $r=\max_{i}\vert\lambda_{i}\vert$ (and of course on $X$ and
$Y$). In this generality, we cannot discuss the question of the sharpness
of the bounds obtained. First, we prove the following lemma.

\begin{lem}\label{lemma311} Under $\left(P_{2}\right)$,
$\left(P_{3}\right)$ {\small and }$\left(P_{4}\right)$, $B_{\sigma}X$
is a closed subspace of X and moreover if $\sigma$ is a finite subset
of $\mathbb{D}$, \[
B_{\sigma}X=\left\{ f\in X:\, f\left(\lambda\right)=0,\,\forall\lambda\in\sigma\,(including\: multiplicities)\right\} .\]
\end{lem}
\begin{proof}
Since $X\subset{\rm Hol}(\mathbb{D})$ continuously, and evaluation
functionals $f\mapsto f(\lambda)$ and $f\mapsto f^{(k)}(\lambda),\, k\in\mathbb{N}^{\star}$,
are continuous on ${\rm Hol}(\mathbb{D}),$ the subspace \[
M=\left\{ f\in X:\, f\left(\lambda\right)=0,\,\forall\lambda\in\sigma\,(including\: multiplicities)\right\} ,\]
 is closed in $X.$

On the other hand, $B_{\sigma}X\subset X,$ and hence $B_{\sigma}X\subset M.$
Indeed, properties $\left(P_{2}\right)$ and $\left(P_{3}\right)$
imply that\textbf{ }$h.X\subset X$, for all $h\in{\rm Hol}((1+\epsilon)\mathbb{D})$\textbf{
}with\textbf{ $\epsilon>0.$} We can write $h=\sum_{k\geq0}\widehat{h}(k)z^{k}$
with $\left|\widehat{h}(k)\right|\leq Cq^{n}$, $C>0$ and $q<1$.
Then $\sum_{n\geq0}\left\Vert \widehat{h}(k)z^{k}f\right\Vert _{X}<\infty$
for every $f\in X.$ Since $X$ is a Banach space we get $hf=\sum_{n\geq0}\widehat{h}(k)z^{k}f\in X$.

\vspace{0.1cm}

In order to see that $M\subset B_{\sigma}X$, it suffices to justify
that \[
\left[f\in X\:{\rm and}\: f(\lambda)=0\right]\Longrightarrow\left[f/b_{\lambda}=(1-\overline{\lambda}z)f/(\lambda-z)\in X\right],\]
 but this is obvious from $\left(P_{4}\right)$ and the previous arguments.
\end{proof}
In Definition \ref{def312} below, $\sigma=\left\{ \lambda_{1},\,...,\,\lambda_{n}\right\} $
is a sequence in the unit disc $\mathbb{D}$ and $B_{\sigma}$ is
the corresponding Blaschke product.
\begin{dfn}\label{def312}\textit{The model space $K_{B_{\sigma}}$. }We define $K_{B_{\sigma}}$ to
be the $n$-dimensional space:
\def\theequation{${3.1.1}$}\begin{equation}
K_{B_{\sigma}}=\left(B_{\sigma}H^{2}\right)^{\perp}=H^{2}\ominus B_{\sigma}H^{2}.\label{eq:3.1.1}\end{equation}
\textit{Malmquist basis}. For $k\in[1,\, n]$, we set $f_{k}(z)=\frac{1}{1-\overline{\lambda_{k}}z},$
and define the family $\left(e_{k}\right)_{1\leq k\leq n}$, (which
is known as \textit{Malmquist basis}, see \cite{Nik3} p.117), by
\def\theequation{${3.1.2}$}\begin{equation}
e_{1}=\frac{f_{1}}{\left\Vert f_{1}\right\Vert _{2}}\,\,\,{\rm and}\,\,\, e_{k}=\left({\displaystyle \prod_{j=1}^{k-1}}b_{\lambda_{j}}\right)\frac{f_{k}}{\left\Vert f_{k}\right\Vert _{2}}\,,\label{eq:3.1.2}\end{equation}
for $k\in[2,\, n]$, where $\left\Vert f_{k}\right\Vert _{2}=\left(1-\vert\lambda_{k}\vert^{2}\right)^{-1/2}.$
The \textit{Malmquist family} $\left(e_{k}\right)_{1\leq k\leq n}$
corresponding to $\sigma$ is an orthonormal basis of $K_{B_{\sigma}}$.

\textit{The orthogonal projection ~$P_{B_{\sigma}}$on $K_{B_{\sigma}}.$
}We define $P_{B_{\sigma}}$ to be the orthogonal projection of $H^{2}$
on its $n$-dimensional subspace $K_{B_{\sigma}}.$ In particular,

\def\theequation{${3.1.3}$}\begin{equation}
P_{B_{\sigma}}=\sum_{k=1}^{n}\left(\cdot,\, e_{k}\right)_{H^{2}}e_{k},\label{eq:3.1.3}\end{equation}
where $\left(\cdot,\,\cdot\right)_{H^{2}}$ means the scalar product
on $H^{2}$. \end{dfn}

\begin{lem}\label{lemma316} Let $\sigma=\left\{ \lambda_{1},\,...,\,\lambda_{n}\right\} $
be a sequence in the unit disc $\mathbb{D}$ , $\left(e_{k}\right)_{1\leq k\leq n}$
be the Malmquist family corresponding to $\sigma$ and $\left\langle \cdot,\,\cdot\right\rangle $
be the Cauchy sesquilinear form $\left\langle h,\, g\right\rangle =\sum_{j\geq0}\hat{h}(j)\overline{\hat{g}(j)}$
for $h,\, g\in{\rm Hol}(\mathbb{D}).$ For all $f\in{\rm Hol}(\mathbb{D})$
and $1\leq k\leq n,$ the series $\left\langle f,\, e_{k}\right\rangle $
are absolutely converging. Moreover, if $Y$ is a Banach space satisfying
$\left(P_{1}\right),$ the map $\Phi:\,{\rm Hol}(\mathbb{D})\rightarrow Y\subset{\rm Hol}(\mathbb{D})$
given by
\[\Phi:\,\, f\mapsto\sum_{k=1}^{n}\left\langle f,\, e_{k}\right\rangle e_{k},\]
is well defined and has the following properties:

(a) \textbf{$\Phi_{\vert H^{2}}=P_{B_{\sigma}},$}

(b) $\Phi$ is continuous on ${\rm Hol}(\mathbb{D})$ with the topology
of the uniform convergence on compact sets of $\mathbb{D}$,

\begin{onehalfspace}
(c) if $X$ satisfies $\left(P_{2}\right)$, $\left(P_{3}\right)$,{\small{}
}$\left(P_{4}\right)$ and $Y\subset X,$ and if $\Psi=Id_{\vert X}-\Phi_{\vert X},$
then ${\rm Im}\left(\Psi\right)\subset B_{\sigma}X,$
\end{onehalfspace}

(d) if $f\in{\rm Hol}(\mathbb{D}),$ then \[
\left|\Phi(f)(\zeta)\right|=\left|\left\langle f,\, P_{B_{\sigma}}k_{\zeta}\right\rangle \right|=\left|\left\langle f,\,\sum_{k=1}^{n}\overline{e_{k}(\zeta)}e_{k}\right\rangle \right|,\]
for all $\zeta\in\mathbb{D},$ where $P_{B_{\sigma}}$ is defined
in (\ref{eq:3.1.2}) and $k_{\zeta}=\left(1-\overline{\zeta}z\right)^{-1}.$\end{lem}
\begin{proof}

First of all we set $r=\max{}_{\lambda\in\sigma}\left|\lambda\right|.$
If $f\in{\rm Hol}(\mathbb{D})$ and $\rho\in]0,1[,$ then\[
\widehat{f}(j)=\left(2\pi\right)^{-1}\int_{\rho\mathbb{T}}f(w)w^{-j-1}dw,\]
 for all $j\geq0.$ For a subset $A$ of $\mathbb{C}$ and for a bounded
function $h$ on $A$, we define $\left\Vert h\right\Vert _{A}:=\sup_{z\in A}\left|h(z)\right|.$
As a result,

\def\theequation{${3.1.4}$}\begin{equation}
\left|\left\langle f,\, e_{k}\right\rangle \right|\leq\sum_{j\geq0}\left|\widehat{f}(j)\overline{\widehat{e_{k}}(j)}\right|\leq\left(2\pi\rho\right)^{-1}\left\Vert f\right\Vert _{\rho\mathbb{T}}\sum_{j\geq0}\left|\widehat{e_{k}}(j)\right|\rho^{-j}.\label{eq:3.1.4}\end{equation}
Now if $\rho$ is close enough to $1$, it satisfies the inequality
$1\leq\rho^{-1}<r^{-1}$, which entails $\sum_{j\geq0}\left|\widehat{e_{k}}(j)\right|\rho^{-j}<+\infty$
for each $k=1,\,...,\, n$, and the series $\left\langle f,\, e_{k}\right\rangle $
are absolutely converging.

Now, the the point (a) is a direct consequence of (\ref{eq:3.1.3}). In order to
check point (b), let $\left(f_{l}\right)_{l\in\mathbb{N}}$ be a sequence
of ${\rm Hol}(\mathbb{D})$ converging to $0$ uniformly on compact
sets of $\mathbb{D}$. We need to see that $\left(\Phi\left(f_{l}\right)\right)_{l\in\mathbb{N}}$
converges to $0$, for which it is sufficient to show that $\lim_{l}\left|\left\langle f_{l},\, e_{k}\right\rangle \right|=0,$
for every $k=1,\,2,\,...,\, n,$ which is clear applying (\ref{eq:3.1.4}) to
$f=f_{l}.$

We now prove point (c). Using point (a), since $Pol_{+}\subset H^{2}$,
we get that ${\rm Im}\left(\Psi_{\vert Pol_{+}}\right)\subset B_{\sigma}H^{2}$.
Now, since $Pol_{+}\subset Y$ and ${\rm Im}(\Phi)\subset Y$, we
deduce that \[
{\rm Im}\left(\Psi_{\vert Pol_{+}}\right)\subset B_{\sigma}H^{2}\cap Y\subset B_{\sigma}H^{2}\cap X,\]
since $Y\subset X$. Now $\Psi\left(p\right)\in X$ and satisfies
$\left(\Psi\left(p\right)\right)_{\vert\sigma}=0$ (that is to say
$\left(\Psi\left(p\right)\right)\left(\lambda\right)=0,\,\forall\lambda\in\sigma$
(${\rm including\: multiplicities}$)\textbf{ }for all $p\in Pol_{+}.$
Using Lemma \ref{lemma311}, we get that ${\rm Im}\left(\Psi_{\vert Pol_{+}}\right)\subset B_{\sigma}X.$
Now, $Pol_{+}$ being dense in $X$ (property $\left(P_{2}\right)$),
and $\Psi$ being continuous on $X$ (point (b)), we can conclude
that ${\rm Im}\left(\Psi\right)\subset B_{\sigma}X$.

In order to prove (d), we simply need to write that \[
\Phi(f)(\zeta)=\sum_{k=1}^{n}\left\langle f,\, e_{k}\right\rangle e_{k}(\zeta)=\left\langle f,\,\sum_{k=1}^{n}\overline{e_{k}(\zeta)}e_{k}\right\rangle ,\]
\textit{$\forall f\in{\rm Hol}(\mathbb{D}),$ $\forall\zeta\in\mathbb{D}$
}and to notice that \textit{$\sum_{k=1}^{n}\overline{e_{k}(\zeta)}e_{k}=\sum_{k=1}^{n}\left(k_{\zeta},\, e_{k}\right)_{H^{2}}e_{k}=P_{B_{\sigma}}k_{\zeta}.$}
\end{proof}

\begin{proof}[Proof of Theorem \ref{Theorem A} ] Let $\sigma=\left\{ \lambda_{1},\,...,\,\lambda_{n}\right\} $
be a sequence in the unit disc $\mathbb{D}$ and $\left(e_{k}\right)_{1\leq k\leq n}$
the Malmquist family (\ref{eq:3.1.2}) associated to $\sigma$. Taking $f\in X$
, we set $g=\Phi(f)=\sum_{k=1}^{n}\left\langle f,\, e_{k}\right\rangle e_{k},$
where $\Phi$ is defined in Lemma \ref{lemma316}. In the same spirit of that
of (\ref{eq:3.1.4}), we notice that
\[
\widehat{e_{k}}(j)=\left(2\pi i\right)^{-1}\int_{R\mathbb{T}}e_{k}(w)w^{-j-1}{\rm d}w,\]
for all $j\geq0$ and for all $R,$ $1<R<\frac{1}{r}$. As a result,
\[
\left|\widehat{e_{k}}(j)\right|\leq\left(2\pi R^{j+1}\right)^{-1}\left\Vert e_{k}\right\Vert _{R\mathbb{T}}\,\,\mbox{and}\]
\[
\,\,\sum_{j\geq0}\left|\hat{f}(j)\overline{\widehat{e_{k}}(j)}\right|
\leq(2\pi R)^{-1}\left\Vert e_{k}\right\Vert _{R\mathbb{T}}\sum_{j\geq0}\left|\hat{f}(j)\right|R^{-j}<\infty,\]
since $R>1$ and $f$ is holomorphic in $\mathbb{D},$ (where $\left\Vert \cdot\right\Vert _{A}$
is defined above in the proof of Lemma \ref{lemma316}).

We now suppose that $\left\Vert f\right\Vert _{X}\leq1$. Since ${\rm Hol}\left(r^{-1}\mathbb{D}\right)\subset Y$,
we have $g\in Y$ and using Lemma \ref{lemma316} point (c) we get
\[
f-g=\Psi(f)\in B_{\sigma}X,\]
where $\Psi$ is defined in Lemma \ref{lemma316}, as $\Phi$. Moreover,
\[
\left\Vert g\right\Vert _{Y}\leq\sum_{k=1}^{n}\left|\left\langle f,e_{k}\right\rangle _{}\right|\left\Vert e_{k}\right\Vert _{Y}.\]
In order to bound the right-hand side, recall that for all $j\geq0$
and for $R=2/(r+1)\in]1,1/r[$,
\[
\sum_{j\geq0}\left|\widehat{f}(j)\overline{\widehat{e_{k}}(j)}\right|\leq(2\pi)^{-1}\left\Vert e_{k}\right\Vert _{2(r+1)^{-1}\mathbb{T}}\sum_{j\geq0}\left|\widehat{f_{}}(j)\right|\left(2^{-1}(r+1)\right)^{j}.\]
Since the norm $f\mapsto\sum_{j\geq0}\left|\widehat{f}(j)\right|\left(2^{-1}(r+1)\right)^{j}$
is continuous on ${\rm Hol}(\mathbb{D}),$ and the inclusion $X\subset{\rm Hol}(\mathbb{D})$
is also continuous, there exists $C_{r}>0$ such that
\[
\sum_{j\geq0}\left|\widehat{f}(j)\right|\left(2^{-1}(r+1)\right)^{j}\leq C_{r}\parallel f\parallel_{X},\]
for every $f\in X$. On the other hand, ${\rm Hol}\left(2(r+1)^{-1}\mathbb{D}\right)\subset Y$
(continuous inclusion again), and hence there exists $K_{r}>0$ such
that
\[
\left\Vert e_{k}\right\Vert _{Y}\leq K_{r}\sup_{\vert z\vert<2(r+1)^{-1}}\left|e_{k}(z)\right|=K_{r}\left\Vert e_{k}\right\Vert _{2(r+1)^{-1}\mathbb{T}}.\]

It is more or less clear that the right-hand side of the last inequality
can be bounded in terms of $r$ and $n$ only. Let us give a proof
to this fact. It is clear that it suffices to estimate
\[
\sup_{1<\vert z\vert<2(r+1)^{-1}}\left|e_{k}(z)\right|.\]
In order to bound this quantity, notice that

\def\theequation{3.1.5}
\begin{equation}
\vert b_{\lambda}(z)\vert^{2}\leq\left|\frac{\lambda-z}{1-\bar{\lambda}z}\right|^{2}=1+\frac{(\vert z\vert^{2}-1)(1-\vert\lambda\vert^{2})}{\vert1-\bar{\lambda}z\vert^{2}},\label{eq:3.1.5}\end{equation}
for all $\lambda\in\mathbb{D}$ and all $z\in\vert\lambda\vert^{-1}\mathbb{D}$.
Using the identity (\ref{eq:3.1.5}) for $\lambda=\lambda_{j}$, $1\leq j\leq n$,
and $z=\rho e^{it}$, $\rho=2(1+r)^{-1}$, we get
\[
\left|e_{k}(\rho e^{it})\right|^{2}\leq\left(\prod_{j=1}^{k-1}\left|b_{\lambda_{j}}(\rho e^{it})\right|^{2}\right)\left|\frac{1}{1-\overline{\lambda_{k}}\rho e^{it}}\right|^{2}\leq\]
\[
\leq\left(\prod_{j=1}^{k-1}\left(1+\frac{(\rho^{2}-1)(1-\vert\lambda_{j}\vert^{2})}{1-\vert\lambda_{j}\vert^{2}\rho^{2}}\right)\right)\left(\frac{1}{1-\vert\lambda_{k}\vert\rho}\right)^{2},\]
for all $k=2,\,...,\, n$. Expressing $\rho$ in terms of $r$, we
obtain
\[
\left\Vert e_{k}\right\Vert _{2(r+1)^{-1}\mathbb{T}}\leq\frac{1}{1-\frac{2r}{r+1}}\sqrt{2\left(\prod_{j=1}^{n-1}\left(1+\frac{2(\frac{1}{r^{2}}-1)}{1-r^{2}\frac{4}{(r+1)^{2}}}\right)\right)}=:C_{1}(r,\, n),\]
and
\[
\sum_{j\geq0}\left|\hat{f}(j)\overline{\hat{e_{k}}(j)}\right|\leq(2\pi)^{-1}C_{r}\left\Vert e_{k}\right\Vert {}_{2(r+1)^{-1}\mathbb{T}}\parallel f\parallel_{X}\leq(2\pi)^{-1}C_{r}C_{1}(r,\, n)\parallel f\parallel_{X}.\]
On the other hand, since
\[
\left\Vert e_{k}\right\Vert _{Y}\leq K_{r}\left\Vert e_{k}\right\Vert _{2(r+1)^{-1}\mathbb{T}}\leq K_{r}C_{1}(r,\, n),\]
we get
\[
\left\Vert g\right\Vert _{Y}\leq\sum_{k=1}^{n}(2\pi)^{-1}C_{r}C_{1}(r,\, n)\left\Vert f\right\Vert _{X}K_{r}C_{1}(r,n)=(2\pi)^{-1}nC_{r}K_{r}\left(C_{1}(r,\, n)\right)^{2}\left\Vert f\right\Vert _{X},\]
which proves that

\texttt{\textbf{\textit{\[
c(\sigma,\, X,\, Y)\leq(2\pi)^{-1}nC_{r}K_{r}\left(C_{1}(r,\, n)\right)^{2},\]
}}}and completes the proof of Theorem \ref{Theorem A}.\end{proof}

\subsection{The case where $X$ is a Hilbert space}\label{subsection32}

We suppose in this subsection that $X$ is a Hilbert space and that
$X,\, Y$ satisfy properties $\left(P_{i}\right)_{1\leq i\leq4}.$
We prove Theorem \ref{Theorem B}  and obtain a better estimate for $c\left(\sigma,\, X,\, Y\right)$
than in Theorem \ref{Theorem A}
 (see point (i) of Theorem \ref{Theorem B} ). For the case $Y=H^{\infty}$,
(point (ii) of Theorem \ref{Theorem B} ), we can considerably improve this estimate.
We omit an easy proof of the following lemma.

\begin{lem}\label{lemma321} Let $\sigma=\{\lambda_{1},...,\,\lambda_{1},\,\lambda_{2},\,...,\,\lambda_{2},...,\,\lambda_{t},...,\,\lambda_{t}\}$
be a finite sequence of $\mathbb{D}$ where every $\lambda_{s}$ is
repeated according to its multiplicity $m_{s}$, $\sum_{s=1}^{t}m_{s}=n$.
Let $\left(H,\,\left(\cdot,\,\cdot\right)_{H}\right)$ be a Hilbert
space continuously embedded into ${\rm Hol}(\mathbb{D})$ and satisfying
properties $\left(P_{i}\right)$ for $i=2,\,3,\,4$. Then \[
K_{B_{\sigma}}(H)=:H\ominus B_{\sigma}H={\rm span}\left(k_{\lambda_{j},\, i}:\,1\leq j\leq t,\,0\leq i\leq m_{j}-1\right),\]
where $k_{\lambda,\, i}=\left(\frac{d}{d\overline{\lambda}}\right)^{i}k_{\lambda}$
and $k_{\lambda}$ is the reproducing kernel of $H$ at point $\lambda$
for every $\lambda\in\mathbb{D},$ i.e. $k_{\lambda}\in H$ and $f\left(\lambda\right)=\left(f,\, k_{\lambda}\right)_{H},\;\forall f\in H.$\end{lem}

\begin{proof}[Proof of Theorem \ref{Theorem B} ]  i). Let $f\in X$, $\left\Vert f\right\Vert _{X}\leq1$
. Lemma \ref{lemma321} shows that
\[
g=P_{B_{\sigma}}^{H}f=\sum_{k=1}^{n}\left(f,\,\mathcal{E}_{k}\right)_{H}\mathcal{E}_{k}\]
is the orthogonal projection of $f$ onto subspace $K_{B_{\sigma}}$.
Function $g$ belongs to $Y$ because all $k_{\lambda_{j},i}$ are
in ${\rm Hol}((1+\epsilon)\mathbb{D})$ for a convenient $\epsilon>0$,
and $Y$ satisfies $(P_{1}).$

On the other hand, $g-f\in B_{\sigma}H$ (again by Lemma \ref{lemma321}). Moreover,
using Cauchy-Schwarz inequality,
\[
\left\Vert g\right\Vert _{Y}\leq\sum_{k=1}^{n}\left|\left(f,\,\mathcal{E}_{k}\right)_{H}\right|\left\Vert \mathcal{E}_{k}\right\Vert _{Y}\leq\left(\sum_{k=1}^{n}\left|\left(f,\,\mathcal{E}_{k}\right)_{H}\right|^{2}\right)^{1/2}\left(\sum_{k=1}^{n}\left\Vert \mathcal{E}_{k}\right\Vert _{Y}^{2}\right)^{1/2}\leq\]
\[
\leq\left\Vert f\right\Vert _{H}\left(\sum_{k=1}^{n}\left\Vert \mathcal{E}_{k}\right\Vert _{Y}^{2}\right)^{1/2},\]
which proves i).

ii). If $Y=H^{\infty}$, then
\[
\vert g(\zeta)\vert=\left|\left(P_{B_{\sigma}}^{H}f,\, k_{\zeta}\right)_{H}\right|=\left|\left(f,\, P_{B_{\sigma}}^{H}k_{\zeta}\right)_{H}\right|\leq\left\Vert f\right\Vert _{H}\left\Vert P_{B_{\sigma}}^{H}k_{\zeta}\right\Vert _{H},\]
for all $\zeta\in\mathbb{D}$, which proves ii).
\end{proof}

\section{Upper bounds for $C_{n,\, r}\left(H^{2},\, H^{\infty}\right)$}\label{section4}

Here, we specialize the upper estimate obtained in point (ii) of Theorem \ref{Theorem B}  for the case $X=H^{2}$, the Hardy space of the disc. Later on,
we will see that this estimate is sharp at least for some special
sequences $\sigma$ (see Section 7). We also develop a slightly different
approach to the interpolation constant $c\left(\sigma,\, H^{2},\, H^{\infty}\right)$
giving more estimates for individual sequences $\sigma=\left\{ \lambda_{1},\,...,\,\lambda_{n}\right\} $
of $\mathbb{D}$. We finally prove the right-hand side inequality
of Theorem \ref{Theorem C} for the particular case $p=2$.

\begin{prop}\label{prop41}For every sequence
$\sigma=\{\lambda_{1},\,...,\,\lambda_{n}\}$ of $\mathbb{D}$ we
have
\def\theequation{$4.1.1$}
\begin{equation}
c\left(\sigma,\, H^{2},\, H^{\infty}\right)\leq\sup_{\zeta\in\mathbb{D}}\left(\frac{1-\vert B_{\sigma}(\zeta)\vert^{2}}{1-\vert\zeta\vert^{2}}\right)^{1/2},\label{eq:4.1.1}\end{equation}
and therefore\def\theequation{$4.1.2$} \end{prop}
\def\theequation{$4.1.2$}
\begin{equation}
c\left(\sigma,\, H^{2},\, H^{\infty}\right)\leq\sqrt{2}\sup_{\vert\zeta\vert=1}\left|B'(\zeta)\right|^{\frac{1}{2}}=\sqrt{2}\sup_{\vert\zeta\vert=1}\left|\sum_{i=1}^{n}\frac{1-\vert\lambda_{i}\vert^{2}}{\left(1-\bar{\lambda_{i}}\zeta\right)^{2}}\frac{B_{\sigma}(\zeta)}{b_{\lambda_{i}}(\zeta)}\right|^{1/2}.\label{eq:4.1.2}
\end{equation}

\begin{proof}

We prove (\ref{eq:4.1.1}). In order to simplify the notation, we set $B=B_{\sigma}$.
Applying point (ii) of Theorem \ref{Theorem B}  for $X=H^{2}$ and $Y=H^{\infty},$
and using \[
k_{\zeta}(z)=\frac{1}{1-\bar{\zeta}z}\,\,\:\mbox{and}\,\,\:\left(P_{B_{\sigma}}k_{\zeta}\right)(z)=\frac{1-\overline{B_{\sigma}(\zeta)}B_{\sigma}(z)}{1-\overline{\zeta}z},\]
 (see \cite{Nik3} p.199), we obtain
\[
\left\Vert P_{B_{\sigma}}k_{\zeta}\right\Vert _{H^{2}}=\left(\frac{1-\vert B_{\sigma}(\zeta)\vert^{2}}{1-\vert\zeta\vert^{2}}\right)^{1/2},\]
which gives the result.

We now prove (\ref{eq:4.1.2}) using (\ref{eq:4.1.1}). The map $$\zeta\mapsto\left\Vert P_{B}\left(k_{\zeta}\right)\right\Vert =\sup\left\{ \left|f(\zeta)\right|:\, f\in K_{B},\,\left\Vert f\right\Vert \leq1\right\},$$
and hence the map
\[
\zeta\mapsto\left(\frac{1-\vert B(\zeta)\vert^{2}}{1-\vert\zeta\vert^{2}}\right)^{1/2},\]
is a subharmonic function so

\begin{onehalfspace}
\textit{\[
\sup_{\vert\zeta\vert<1}\left(\frac{1-\vert B(\zeta)\vert^{2}}{1-\vert\zeta\vert^{2}}\right)^{1/2}\leq\sup_{\vert w\vert=1}\lim_{r\rightarrow1}\left(\frac{1-\vert B(rw)\vert^{2}}{1-\vert rw\vert^{2}}\right)^{1/2}.\]
}Now applying Taylor's Formula of order 1 at points $w\in\mathbb{T}$
and $u=rw$, $0<r<1$ (it is applicable because $B$ is holomorphic
at every point of $\mathbb{T}$), we get
\end{onehalfspace}
\[
\left(B(u)-B(w)\right)(u-w)^{-1}=B'(w)+o(1),\]
and since $\vert u-w\vert=1-\vert u\vert,$
\[
\left|\left(B(u)-B(w)\right)(u-w)^{-1}\right|=\vert B(u)-B(w)\vert\left(1-\vert u\vert\right)^{-1}=\vert B'(w)+o(1)\vert.\]
Then we have
\[
\vert B(u)-B(w)\vert\geq\vert B(w)\vert-\vert B(u)\vert=1-\vert B(u)\vert,\]
\[
\left(1-\vert B(u)\vert\right)\left(1-\vert u\vert\right)^{-1}\leq\left(1-\vert u\vert\right)^{-1}\vert B(u)-B(w)\vert=\vert B'(w)+o(1)\vert,\]
and
\[
\lim_{r\rightarrow1}\left(\left(1-\vert B(rw)\vert\right)(1-\vert rw\vert)^{-1}\right)^{\frac{1}{2}}\leq\sqrt{\vert B'(w)\vert}.\]
Moreover,
\[
B'(w)=-\sum_{i=1}^{n}\left(1-\vert\lambda_{i}\vert^{2}\right)(1-\overline{\lambda_{i}}w)^{-2}\prod_{j=1,\, j\neq i}^{n}b_{\lambda_{j}}(w),\]
for all $w\in\mathbb{T}$ . This completes the proof since
\[
\frac{1-\vert B(rw)\vert^{2}}{1-\vert rw\vert^{2}}=\frac{(1-\vert B(rw)\vert)(1+\vert B(rw)\vert)}{(1-\vert rw\vert)(1+\vert rw\vert)}\leq2\frac{1-\vert B(rw)\vert}{1-\vert rw\vert}.\]
\end{proof}
\begin{cor}\label{cor42} Let $n\geq1$ and
$r\in[0,1[$. Then,
\[
C_{n,\, r}(H^{2},H^{\infty})\leq2\left(n(1-r)^{-1}\right)^{\frac{1}{2}}.\]
 \end{cor}
\begin{proof}
Indeed, applying Proposition \ref{prop41}
 we obtain \[
\left|B'(w)\right|\leq\left|\sum_{i=1}^{n}\frac{1-\vert\lambda_{i}\vert^{2}}{\left(1-\vert\lambda_{i}\vert\right)^{2}}\right|\leq n\frac{1+r}{1-r}\leq\frac{2n}{1-r}.\]
\end{proof}

Now, we develop a slightly different approach to the interpolation
constant $c\left(\sigma,\, H^{2},\, H^{\infty}\right)$.

\begin{cor}\label{cor43} For every sequence
$\sigma=\left\{ \lambda_{1},\,...,\,\lambda_{n}\right\} $ of $\mathbb{D}$,
\[
c\left(\sigma,\, H^{2},\, H^{\infty}\right)\leq\sup_{z\in\mathbb{T}}\left(\sum_{k=1}^{n}\frac{\left(1-\vert\lambda_{k}\vert^{2}\right)}{\vert z-\lambda_{k}\vert^{2}}\right)^{1/2}.\]
\end{cor}
\begin{proof}
In order to simplify the notation, we set $B=B_{\sigma}$. We consider
$K_{B}$ (see Definition \ref{def312}) and the Malmquist family $\left(e_{k}\right)_{1\leq k\leq n}$
corresponding to $\sigma$ (see Definition \ref{def312}). Now, let $f\in H^{2}$
and
\[
g=P_{B}f=\sum_{k=1}^{n}\left(f,e_{k}\right)_{H^{2}}e_{k},\]
(see (\ref{eq:3.1.3})). The function $g$ belongs to $H^{\infty}$ (it is a
finite sum of $H^{\infty}$ functions) and satisfies $g-f\in BH^{2}$.
Applying Cauchy-Schwarz inequality we get
\[
\vert g(\zeta)\vert\leq\sum_{k=1}^{n}\left|\left(f,e_{k}\right)_{H^{2}}\right|\left|e_{k}(\zeta)\right|\leq\left(\sum_{k=1}^{n}\left|\left(f,e_{k}\right)_{H^{2}}\right|^{2}\right)^{1/2}\left(\sum_{k=1}^{n}\frac{\left(1-\vert\lambda_{k}\vert^{2}\right)}{\left|1-\lambda_{k}\zeta\right|^{2}}\right)^{1/2},\]
for all $\zeta\in\mathbb{D}$. As a result, since $f$ is an arbitrary
$H^{2}$ function, we obtain
\[
c(\sigma,\, H^{2},\, H^{\infty})\leq\sup_{\zeta\in\mathbb{T}}\left(\sum_{k=1}^{n}\frac{\left(1-\vert\lambda_{k}\vert^{2}\right)}{\vert\zeta-\lambda_{k}\vert^{2}}\right)^{1/2},\]
which completes the proof.
\end{proof}
\begin{cor}\label{cor44}For any sequence $\sigma=\left\{ \lambda_{1},\,...,\,\lambda_{n}\right\} $
in $\mathbb{D}$ ,
\[
c(\sigma,\, H^{2},\, H^{\infty})\leq\left(\sum_{j=1}^{n}\frac{1+\left|\lambda_{j}\right|}{1-\left|\lambda_{j}\right|}\right)^{1/2}.\]
\end{cor}
\begin{proof}
Indeed,
\begin{onehalfspace}
\[
\sum_{k=1}^{n}\frac{\left(1-\vert\lambda_{k}\vert^{2}\right)}{\vert\zeta-\lambda_{k}\vert^{2}}\leq\left(\sum_{k=1}^{n}\frac{\left(1-\vert\lambda_{k}\vert^{2}\right)}{\left(1-\vert\lambda_{k}\vert\right){}^{2}}\right)^{1/2}\]
and the result follows using Theorem 4.3.\end{onehalfspace}
\end{proof}
Now we prove the right-hand side inequality of Theorem \ref{Theorem C}  for the special
case $p=2$.

\begin{proof}[Proof of Theorem \ref{Theorem C}  ($p=2$, the right-hand side inequality only )]
Since $1+\vert\lambda_{j}\vert\leq2$ and $1-\vert\lambda_{j}\vert\geq1-r$
for all $j\in[1,\, n]$, applying Corollary \ref{cor44}  we get
\[
C_{n,r}(H^{2},H^{\infty})\leq\sqrt{2}n^{1/2}(1-r)^{-1/2}.\]
\end{proof}

\begin{rmk}
As a result, we get once more the same estimate for $C_{n,r}(H^{2},\, H^{\infty})$
as in Corollary \ref{cor42} , with the constant $\sqrt{2}$ instead of $2$.
\end{rmk}
It is natural to wonder if it is possible to improve the bound $\sqrt{2}n^{1/2}(1-r)^{-1/2}$.
We return to this question in Section 7 below.

\section{Upper bounds for $C_{n,\, r}\left(H^{p},\, H^{\infty}\right),\, p\geq1$}
\label{section5}

In this section we extend Corollary \ref{cor42}  to all Hardy spaces $H^{p}$:
we prove the right-hand side inequality of Theorem \ref{Theorem C} , $p\neq2$. We
first prove the following lemma.

\begin{lem}\label{lemma51}\textbf{ }Let $n\geq1$
and $0\leq r<1.$ Then,
\[
C_{n,r}(H^{1},H^{\infty})\leq2n(1-r)^{-1}.\]
\end{lem}
\begin{proof}
Let $\sigma$ be a finite subset of $\mathbb{D},$ $\left(e_{k}\right)_{1\leq k\leq n}$
be the Malmquist basis corresponding to $\sigma$ (see Definition \ref{def312}), and
$f\in H^{1}$ such that $\left\Vert f\right\Vert _{H^{1}}\leq1$.
Let also $g=\Phi(f)$ where $\Phi$ is defined in Lemma \ref{lemma316}. Applying
point (d) of Lemma \ref{lemma316}, we get
\[
\vert g(\zeta)\vert\leq\left\Vert f\right\Vert _{H^{1}}\left\Vert \sum_{k=1}^{n}e_{k}\overline{e_{k}(\zeta)}\right\Vert _{H^{\infty}}\leq\left\Vert \sum_{k=1}^{n}e_{k}\overline{e_{k}(\zeta)}\right\Vert _{H^{\infty}}.\]
Since Blaschke factors have modulus 1 on the unit circle, \[
\left\Vert e_{k}\right\Vert _{H^{\infty}}\leq\left(1+\left|\lambda_{k}\right|\right)^{1/2}\left(1-\left|\lambda_{k}\right|\right)^{-1/2}.\]
As a consequence,
{\small
\[
\vert g(\zeta)\vert\leq\sum_{k=1}^{n}\left\Vert e_{k}\right\Vert _{H^{\infty}}\left|\overline{e_{k}(\zeta)}\right|\leq\sum_{k=1}^{n}\left\Vert e_{k}\right\Vert _{H^{\infty}}^{2}\leq\sum_{k=1}^{n}\left(1+\left|\lambda_{k}\right|\right)\left(1-\left|\lambda_{k}\right|\right)^{-1}\leq2n(1-r)^{-1},\]}
for all $\zeta\in\mathbb{D},$ which completes the proof.
\end{proof}
Before proving the upper bound in Theorem \ref{Theorem C} , we give the following
general observation which is a direct consequence of the classical
complex interpolation between Banach spaces, see \cite{Bergh, Tr}. In particular,
we use the notation of \cite{Bergh}, Chapter 4.

\begin{lem}\label{lemma52}\textbf{ }Let $X_{1}$
and $X_{2}$ be two Banach spaces of holomorphic functions in the
unit disc $\mathbb{D}$. Let also $\theta\in[0,\,1]$ and $\left(X_{1},\, X_{2}\right)_{[\theta]}$
be the corresponding intermediate Banach space resulting from the
classical complex interpolation method applied between $X_{1}$ and
$X_{2}$. Then,\[
C_{n,\, r}\left(\left(X_{1},\, X_{2}\right)_{[\theta]},\, H^{\infty}\right)\leq C_{n,\, r}\left(X_{1},\, H^{\infty}\right)^{1-\theta}C_{n,\, r}\left(X_{2},\, H^{\infty}\right)^{\theta},\]
for all $n\geq1,\, r\in[0,\,1).$
\end{lem}

\begin{proof}
Let $X$ be a Banach space of holomorphic functions in the unit disc
$\mathbb{D}$ and let $\sigma=\{\lambda_{1},\,\lambda_{2},\,...,\,\lambda_{n}\}\subset\mathbb{D}$
be a finite subset of the disc. Let $T\::\, X\longrightarrow H^{\infty}/B_{\sigma}H^{\infty}$
be the restriction map defined by \[
Tf=\left\{ g\in H^{\infty}:\: f-g\in B_{\sigma}X\right\} ,\]
for every $f\in X$. Then,
\[
\left\Vert T\right\Vert _{X\rightarrow H^{\infty}/B_{\sigma}H^{\infty}}=c\left(\sigma,\, X,\, H^{\infty}\right).\]
Now, since $\left(X_{1},\, X_{2}\right)_{[\theta]}$ is an exact interpolation
space of exponent $\theta$ (see \cite{Bergh} or \cite{Tr} Theorem 1.9.3-(a),
p.59), we can complete the proof.
\end{proof}
Now we prove the right-hand side inequality of Theorem \ref{Theorem C}  for the remaining
case $p\neq2$.

\begin{proof}[Proof of Theorem \ref{Theorem C}  ($p\neq2$, the right-hand side inequality only )]
In the case $X=H^{p},$ there exists $0\leq\theta\leq1$ such that
$1/p=1-\theta$, and $\left[H^{1},H^{\infty}\right]_{\theta}=H^{p}$
(a topological identity: the spaces are the same and the norms are
equivalent (up to constants depending on $p$ only), see \cite{Jones}).
As a consequence, applying Lemma 5.2 with $X_{1}=H^{1}$ and $X_{2}=H^{\infty},$
we get
\[
C_{n,\, r}\left(H^{p},\, H^{\infty}\right)\leq\gamma_{p}\left(C_{n,\, r}\left(H^{1},H^{\infty}\right)\right)^{1-\theta}\left(C_{n,\, r}\left(H^{\infty},H^{\infty}\right)\right)^{\theta}\leq\]
\[
\leq\gamma_{p}\left(2n(1-r)^{-1}\right)^{1-\theta}\left(1\right)^{\theta},\]
where $\theta=1-1/p,$ and $\gamma_{p}$ is a constant depending only
on $p$. Using Lemma \ref{lemma51}  and the fact that $C_{n,\, r}\left(H^{\infty},\, H^{\infty}\right)\leq1,$
we get\[
C_{n,\, r}\left(H^{p},\, H^{\infty}\right)\leq\gamma_{p}\left(2n(1-r)^{-1}\right)^{1/p},\]
which completes the proof.
\end{proof}

\section{Upper bounds for $C_{n,\, r}\left(L_{a}^{2},\, H^{\infty}\right)$}\label{section6}

Our aim here is to generalize Corollary \ref{cor42}  to the case $X=l_{a}^{2}\left((k+1)^{\alpha}\right)$,
$\alpha\in[-1,\,0]$, the Hardy weighted spaces of all $f=\sum_{k\geq0}\hat{f}(k)z^{k}$
satisfying \[
\left\Vert f\right\Vert _{X}^{2}:=\sum_{k\geq0}\left|\hat{f}(k)\right|^{2}(k+1)^{2\alpha}<\infty.\]
Note that $H^{2}=l_{a}^{2}(1)$ and $L_{a}^{2}=\, l_{a}^{2}\left(\left(k+1\right)^{-\frac{1}{2}}\right).$
We prove the upper bound of Theorem \ref{Theorem D}. The main technical tool used
here is a Bernstein-type inequality for rational functions.

\subsection{Bernstein-type inequalities for rational functions}\label{subsection61}

Bernstein-type inequalities for rational functions were the subject
of a number of papers and monographs (see, for instance, {[}2-3, 7-8,
13{]}). Perhaps, the stronger and closer to ours (Proposition \ref{prop611} )
of all known results are due to K.Dyakonov \cite{Dya1, Dya2}. First, we
prove Proposition \ref{prop611}  below, which tells that if $\sigma=\left\{ \lambda_{1},\,...,\,\lambda_{n}\right\} \subset\mathbb{D}$,
$r={\displaystyle \max_{j}}\left|\lambda_{j}\right|,$ $B=B_{\sigma}$
and $f\in K_{B},$ then
\def\theequation{${\star}$}\begin{equation}
\left\Vert f'\right\Vert _{H^{2}}\leq\alpha_{n,\, r}\left\Vert f\right\Vert _{H^{2}},\label{eq:-1}\end{equation}
where $\alpha_{n,\, r}$ is a constant (explicitly given in Proposition \ref{prop611} ) depending on $n$ and $r$ only such that $0<\alpha_{n,\, r}\leq3\frac{n}{1-r}$. Proposition \ref{prop611}  is in fact a special case ($p=2$) of a K. Dyakonov's
result \cite{Dya1} for $p=2$. It is important to recall that this result
is in its turn, a generalization of M. Levin's inequality {[}13{]}
corresponding to the case $p=\infty$. More precisely, it is proved
in \cite{Dya1} that the norm $\left\Vert D\right\Vert _{K_{B}^{p}\rightarrow H^{p}}$
of the differentiation operator $Df=f'$ on the star-invariant subspace
of the Hardy space $H^{p}$, $K_{B}^{p}:=H^{p}\cap B\overline{zH^{p}}$,
(where the bar denotes complex conjugation) satisfies the following
inequalities \[
c'_{p}\left\Vert B'\right\Vert _{\infty}\leq\left\Vert D\right\Vert _{K_{B}^{p}\rightarrow H^{p}}\leq c_{p}\left\Vert B'\right\Vert _{\infty},\]
for every $p$, $1\leq p\leq\infty$ where $c_{p}$ and $c'_{p}$
are positives constants depending on $p$ only, $B$ is a finite Blaschke
product and $\left\Vert \cdot\right\Vert _{\infty}$ means the norm
in $L^{\infty}(\mathbb{T})$. For the partial case considered in Proposition \ref{prop611}  below, our proof is different from \cite{Dya1, Dya2}. It is based
on an elementary Hilbert space construction for an orthonormal basis
in $K_{B}$ and the constant $c_{2}$ obtained is slightly better.
More precisely, it is proved in \cite{Dya1} that $c'_{2}=\frac{1}{36c}$,
$c_{2}=\frac{36+c}{2\pi}$ and $c=2\sqrt{3\pi}$ (as one can check
easily ($c$ is not evaluated in \cite{Dya1})). It implies an inequality
of type $(\star)$ (with a constant about $6.5$ instead of $3$).

In \cite{Z1}, we discuss the {}``asymptotic sharpness'' of our constant
$\alpha_{n,\, r}.$ We find an inequality for \[
C_{n,\, r}=\sup\,\left\Vert D\right\Vert _{K_{B}\rightarrow H^{2}},\]
(sup is over all $B$ with given $n={\rm deg}\, B$ and $r=\max_{\lambda\in\sigma}\left|\lambda\right|$),
which is asymptotically sharp as $n\rightarrow\infty$. Our result
in {[}18{]} is that $\lim_{n\rightarrow\infty}\frac{C_{n,\, r}}{n}=\frac{1+r}{1-r}$
for every $r,\;0\leq r<1$.

\begin{prop}\label{prop611} Let $B=\prod_{j=1}^{n}b_{\lambda_{j}}$,
be a finite Blaschke product (of order n), $r={\displaystyle \max_{j}}\left|\lambda_{j}\right|,$
and $f\in K_{B}=H^{2}\ominus BH^{2}.$ Then for every $n\geq2$ and
$r\in[0,\,1),$\[
\left\Vert f'\right\Vert _{H^{2}}\leq\alpha_{n,\, r}\left\Vert f\right\Vert _{H^{2}},\]
where $\alpha_{n,\, r}=\left[1+(1+r)(n-1)+\sqrt{n-2}\right](1-r)^{-1},$
and in particular\[
\left\Vert f'\right\Vert _{H^{2}}\leq3\frac{n}{1-r}\left\Vert f\right\Vert _{H^{2}},\]
for all $n\geq1$ and $r\in[0,\,1).$ \end{prop}
\begin{proof}
Using (\ref{eq:3.1.3}), $f=P_{B}f=\sum_{k=1}^{n}\left(f,\, e_{k}\right)_{H^{2}}e_{k}$,
for all $f\in K_{B}.$ Noticing that, $e'_{k}$ \[
e_{k}'=\sum_{i=1}^{k-1}\frac{b_{\lambda_{i}}^{'}}{b_{\lambda_{i}}}e_{k}+\overline{\lambda_{k}}\frac{1}{\left(1-\overline{\lambda_{k}}z\right)}e_{k},\]
for $k\in[2,\, n],$ we get\[
f'=\sum_{k=2}^{n}\left(f,\, e_{k}\right)_{H^{2}}\sum_{i=1}^{k-1}\frac{b_{\lambda_{i}}^{'}}{b_{\lambda_{i}}}e_{k}+\sum_{k=1}^{n}\left(f,\, e_{k}\right)_{H^{2}}\overline{\lambda_{k}}\frac{1}{\left(1-\overline{\lambda_{k}}z\right)}e_{k},\]
\[
=\sum_{i=1}^{n}\frac{b_{\lambda_{i}}'}{b_{\lambda_{i}}}\sum_{k=i+1}^{n-1}\left(f,\, e_{k}\right)_{H^{2}}e_{k}+\sum_{k=1}^{n}\left(f,\, e_{k}\right)_{H^{2}}\overline{\lambda_{k}}\frac{1}{\left(1-\overline{\lambda_{k}}z\right)}e_{k}.\]
Now using Cauchy-Schwarz inequality and the fact that $e_{k}$ is
a vector of norm 1 in $H^{2}$ for $k=1,\,...,\, n$, we get \[
\left\Vert \sum_{k=2}^{n}\overline{\lambda_{k}}\left(f,\, e_{k}\right)_{H^{2}}\frac{1}{\left(1-\overline{\lambda_{k}}z\right)}e_{k}\right\Vert _{H^{2}}\leq\sum_{k=1}^{n}\left|\left(f,\, e_{k}\right)_{H^{2}}\right|\left\Vert \overline{\lambda_{k}}\frac{1}{\left(1-\overline{\lambda_{k}}z\right)}\right\Vert _{\infty}\left\Vert e_{k}\right\Vert _{H^{2}}\leq\]
\[
\leq\frac{1}{1-r}\sum_{k=1}^{n}\left|\left(f,\, e_{k}\right)_{H^{2}}\right|\leq\frac{1}{1-r}\left(\sum_{k=1}^{n}\left|\left(f,\, e_{k}\right)_{H^{2}}\right|^{2}\right)^{\frac{1}{2}}\sqrt{n}\leq\frac{1}{1-r}\left\Vert f\right\Vert _{H^{2}}\sqrt{n}.\]
Further,\[
\left\Vert \sum_{i=1}^{n-1}\frac{b_{\lambda_{i}}'}{b_{\lambda_{i}}}\sum_{k=i+1}^{n}e_{k}\left(f,\, e_{k}\right)_{H^{2}}\right\Vert _{H^{2}}\leq\sum_{i=1}^{n-1}\left\Vert \frac{b_{\lambda_{i}}'}{b_{\lambda_{i}}}\right\Vert _{\infty}\left\Vert \sum_{k=i+1}^{n}\left(f,\, e_{k}\right)_{H^{2}}e_{k}\right\Vert _{H^{2}}=\]
\[
=\left(\max_{1\leq i\leq n-1}\left\Vert \frac{b_{\lambda_{i}}'}{b_{\lambda_{i}}}\right\Vert _{\infty}\right)\sum_{i=1}^{n-1}\left(\sum_{k=i+1}^{n}\left|\left(f,\, e_{k}\right)_{H^{2}}\right|^{2}\right)^{\frac{1}{2}}\leq\max_{i}\left\Vert \frac{b_{\lambda_{i}}'}{b_{\lambda_{i}}}\right\Vert _{\infty}\sum_{i=1}^{n-1}\left\Vert f\right\Vert _{H^{2}}.\]
Now, using\[
\left\Vert \frac{b_{\lambda_{i}}'}{b_{\lambda_{i}}}\right\Vert _{\infty}=\left\Vert \frac{\vert\lambda_{i}\vert^{2}-1}{\left(1-\overline{\lambda_{i}}z\right)\left(\lambda_{i}-z\right)}\right\Vert _{\infty}\leq\frac{1+\left|\lambda_{i}\right|}{1-\left|\lambda_{i}\right|}\leq\frac{1+r}{1-r},\]
we get\[
\left\Vert \sum_{i=1}^{n-1}\frac{b_{\lambda_{i}}'}{b_{\lambda_{i}}}\sum_{k=i+1}^{n}\left(f,\, e_{k}\right)_{H^{2}}e_{k}\right\Vert _{H^{2}}\leq(1+r)\frac{n-1}{1-r}\left\Vert f\right\Vert _{H^{2}}.\]
Finally,
\[
\left\Vert f'\right\Vert _{H^{2}}\leq\left[1+(1+r)(n-1)+\sqrt{n-2}\right](1-r)^{-1}\left\Vert f\right\Vert _{H^{2}}.\]
In particular,\[
\left\Vert f'\right\Vert _{H^{2}}\leq\left(2n-2+\sqrt{n}\right)(1-r)^{-1}\left\Vert f\right\Vert _{H^{2}}\leq3n(1-r)^{-1}\left\Vert f\right\Vert _{H^{2}},\]
for all $n\ge1$ and for every $f\in K_{B}.$
\end{proof}

\subsection{An upper bound for $c\left(\sigma,\, L_{a}^{2},\, H^{\infty}\right)$ } \label{subsection62}

Here, we apply Proposition \ref{prop611}  to prove the right-hand side inequality
of Theorem \ref{Theorem D}. We first prove the following corollary.

\begin{cor}\label{cor621} Let $\sigma$
be a sequence in $\mathbb{D}$\textit{.} Then,
\[
c\left(\sigma,\, l_{a}^{2}\left((k+1)^{-1}\right),\, H^{\infty}\right)\leq2\sqrt{10}\left(n(1-r)^{-1}\right)^{3/2}.\]
\end{cor}
\begin{proof}
Indeed, let $X=l_{a}^{2}\left(-1\right)$, $\sigma$ a finite subset
of $\mathbb{D}$ and $B=B_{\sigma}.$ If $f\in X,$ then using part
(c) of Lemma \ref{lemma316}, $f_{|\sigma}$we get that $\Phi(f)_{|\sigma}=f_{|\sigma}.$
Now, denoting $X^{\star}$ the dual of $X$ with respect to the Cauchy
pairing $\left\langle \cdot,\,\cdot\right\rangle $ (defined in Lemma
\ref{lemma316}) and applying point (d) of the same lemma, we obtain $X^{\star}=l_{a}^{2}\left(1\right)$
and\[
\left|\Phi(f)(\zeta)\right|\leq\left\Vert f\right\Vert _{X}\left\Vert P_{B}k_{\zeta}\right\Vert _{X^{\star}}\leq\left\Vert f\right\Vert _{X}K\left(\left\Vert P_{B}k_{\zeta}\right\Vert _{H^{2}}^{2}+\left\Vert \left(P_{B}k_{\zeta}\right)^{'}\right\Vert _{H^{2}}^{2}\right)^{\frac{1}{2}},\]
where \[
K=\max\left\{ 1,\,\sup_{k\geq1}(k+1)k^{-1}\right\} =2.\]
Since $P_{B}k_{\zeta}\in K_{B}$, Proposition \ref{prop611}  implies\[
\left|\Phi(f)(\zeta)\right|\leq\left\Vert f\right\Vert _{X}K\left(\left\Vert P_{B}k_{\zeta}\right\Vert _{H^{2}}^{2}+9\left(n(1-r)^{-1}\right)^{2}\left\Vert P_{B}k_{\zeta}\right\Vert _{H^{2}}^{2}\right)^{\frac{1}{2}}\leq\]
\[\leq A\left(n(1-r)^{-1}\right)^{3/2}\left\Vert f\right\Vert _{X},\]
where $A=K\left(2/2+9\right)^{\frac{1}{2}}=2\sqrt{10},$ since $\left\Vert P_{B}k_{\zeta}\right\Vert _{2}\leq\sqrt{2}\left(n(1-r)^{-1}\right)^{1/2},$
and since we can suppose $n\geq2,$ (the case $n=1$ being obvious).
\end{proof}
Now let us give the proof of the right-hand side inequality of Theorem \ref{Theorem D}.

\begin{proof}[Proof of Theorem \ref{Theorem D} (the right-hand side inequality only)]
The case $\alpha=0$ corresponds to $X=H^{2}$
and has already been studied in Section \ref{section3} (we can choose $A(0)=\sqrt{2}$).
We now suppose $\alpha\in[-1,\,0)$ and we set $\theta=-\alpha,$
$0\leq\theta\leq1,$ so that $\alpha=(1-\theta)(1-N)+\theta.(-N)$.
Since \[
\left(l_{a}^{2}\left(0\right),\, l_{a}^{2}\left(-1\right)\right)_{[\theta]}=l_{a}^{2}\left(\alpha\right),\]
(see \cite{Bergh, Tr}), this gives, using Lemma \ref{lemma52}  with $X_{1}=l_{a}^{2}\left(0\right)=H^{2}$
and $X_{2}=l_{a}^{2}\left(-1\right),$ Corollary \ref{cor621} , that
\[
C_{n,\, r}\left(l_{a}^{2}\left(\alpha\right),\, H^{\infty}\right)\leq\left(\sqrt{2}\left(n(1-r)^{-1}\right)^{\frac{1}{2}}\right)^{1-\theta}\left(2\sqrt{10}\left(n(1-r)^{-1}\right)^{\frac{3}{2}}\right)^{\theta}=\]
\[
=A(0)^{1-\theta}A(1)^{\theta}\left(n(1-r)^{-1}\right)^{\frac{1-\theta}{2}+\frac{3\theta}{2}}.\]
It remains to use $\theta=-\alpha$, set $A(\alpha)=A(0)^{1-\theta}A(1)^{\theta}$
and apply (5.1.2) with $X=l_{a}^{2}\left((k+1)^{\alpha}\right)$.
In particular, for $\alpha=-1/2$ we get $(1-\theta)/2+3\theta/2=1$
and \[
A\left(-1/2\right)=A(0)^{(1-1/2)}A(1)^{1/2}=\sqrt{2}^{1/2}(2\sqrt{10})^{1/2}=\sqrt{2}.10^{\frac{1}{4}}.\]\end{proof}

\section{Lower bounds for $C_{n,\, r}(X,\, H^{\infty})$}\label{section7}

\subsection{The cases $X=$ $H^{2}$ and $X=L_{a}^{2}$}\label{subsection71}

Here, we consider the standard Hardy and Begman spaces, $X=H^{2}=l_{a}^{2}(1)$
and $X=L_{a}^{2}=l_{a}^{2}((k+1)^{-1/2})$, where the spaces $l_{a}^{2}\left((k+1)^{\alpha}\right)$
are defined in Section \ref{section6}, and the problem of lower estimates for the
one-point special case $\sigma_{n,\,\lambda}=\underbrace{\{\lambda,\lambda,...,\lambda\}}_{n}$,
$\lambda\in\mathbb{D}$. Recall the definition of our constrained
interpolation constant for this case
\[
c(\sigma_{n,\,\lambda},\, H,\, H^{\infty})=\sup\left\{ \Vert f\Vert_{H^{\infty}/b_{\lambda}^{n}H^{\infty}}:\, f\in H,\,\Vert f\Vert_{H}\leq1\right\} ,\]
where $\Vert f\Vert_{H^{\infty}/b_{\lambda}^{n}H^{\infty}}=\inf\left\{ \Vert f+b_{\lambda}^{n}g\Vert_{\infty}:\, g\in H\right\} $.
Our goal in this subsection is to prove the sharpness of the upper
estimate stated in Theorem \ref{Theorem C}  $(p=2)$ and in Theorem \ref{Theorem D} for the quantities
$C_{n,\, r}\left(H^{2},\, H^{\infty}\right)$ and $C_{n,\, r}\left(L_{a}^{2},\, H^{\infty}\right),$
that is to say, to get the lower bounds of Theorem \ref{Theorem C}  $(p=2)$ and
Theorem \ref{Theorem D}.

In the proof, we use properties of reproducing kernel Hilbert space
on the disc $\mathbb{D}$, see for example \cite{Nik2}. Let us recall
some of them adapting the general setting to special cases $X=l_{a}^{2}\left((k+1)^{\alpha}\right)$.
The reproducing kernel of $l_{a}^{2}\left((k+1)^{\alpha}\right)$,
by definition, is a $l_{a}^{2}\left((k+1)^{\alpha}\right)$-valued
function $\lambda\longmapsto k_{\lambda}^{\alpha}$, $\lambda\in\mathbb{D}$,
such that $\left(f,\, k_{\lambda}^{w}\right)=f(\lambda)$ for every
$f\in l_{a}^{2}\left((k+1)^{-\alpha}\right)$, where $\left(\cdot,\,\cdot\right)$
means the scalar product $\left(h,\, g\right)=\sum_{k\geq0}\hat{h}(k)\overline{\hat{g}(k)}(k+1)^{-2\alpha}.$
Since one has $f(\lambda)=\sum_{k\geq0}\hat{f}(k)\lambda^{k}(k+1)^{2\alpha}(k+1)^{-2\alpha}$
($\lambda\in\mathbb{D}$), it follows that
\[
k_{\lambda}^{\alpha}(z)={\displaystyle \sum_{k\geq0}(k+1)^{2\alpha}{\displaystyle \overline{\lambda}^{k}z^{k}}},\: z\in\mathbb{D}.\]
In particular, for the Hardy space $H^{2}=l_{a}^{2}(1)$ ($\alpha=0$),
we get the Szegö kernel $k_{\lambda}(z)=(1-\overline{\lambda}z)^{-1}$
and for the Bergman space $L_{a}^{2}=l_{a}^{2}\left(\left(k+1\right)^{-1/2}\right)$
($\alpha=-1/2$), the Bergman kernel $k_{\lambda}^{-1/2}(z)=$ $(1-\overline{\lambda}z)^{-2}$.

We will use the previous observations for the following composed reproducing
kernels (Aronszajn-deBranges, see \cite{Aro, Nik2}): given the reproducing
kernel $k$ of $H^{2}$ and $\varphi\in\{z^{N}:\: N=1,\,2\},$ the
function $\varphi\circ k$ is also positive definite and the corresponding
Hilbert space is
\[
H_{\varphi}=\varphi(H^{2})=l_{a}^{2}\left((k+1)^{\frac{1-N}{2}}\right).\]
It satisfies the following property: for every $f\in H^{2}$, $\varphi\circ f\in\varphi(H^{2})$
and $\Vert\varphi\circ f\Vert_{\varphi(H^{2})}^{2}\leq\varphi(\Vert f\Vert_{H^{2}}^{2})$
(see \cite{Nik2}, page 320).

We notice in particular that

\def\theequation{${7.1.1}$}\begin{equation}
H_{z}=H^{2}\:\:\mbox{and}\:\: H_{z^{2}}=L_{a}^{2}.\label{eq:-1}\end{equation}
The above relation between the spaces $H^{2}$, $L_{a}^{2}$ and the
spaces $\varphi(H^{2})=H_{\varphi}$ leads to establish the proof
of the left-hand side inequalities stated in Theorem \ref{Theorem C}  (for $p=2$
only) and in Theorem \ref{Theorem D} .
\begin{proof}[Proof of Theorem \ref{Theorem C} ($p=2$) and Theorem \ref{Theorem D}](left-hand side inequalities only)

0)
For $N=1,\,2$ we set  $\varphi_{N}(z)=z^{N}.$

\vspace{0.2cm}
1) We set\[
Q_{n}=\sum_{k=0}^{n-1}(1-\vert\lambda\vert^{2})^{1/2}b_{\lambda}^{k}\left(1-\overline{\lambda}z\right)^{-1},\: H_{n,\, N}=\,\varphi_{N}\circ Q_{n}\,\,\mbox{\,\ and}\,\,\,\Psi=bH_{n,\, N}\,,\; b>0.\]

Then $\Vert Q_{n}\Vert_{2}^{2}=\, n$, and hence by the above Aronszajn-deBranges
inequality,
\[
\Vert\Psi\Vert_{H_{\varphi}}^{2}\leq b^{2}\varphi_{N}\left(\Vert Q_{n}\Vert_{2}^{2}\right)=\, b^{2}\varphi_{N}(n).\]
Let $b>0$ such that $b^{2}\varphi_{N}(n)=\,1$.

\vspace{0.2cm}
 2) Since the spaces $H_{\varphi_{N}}$ and $H^{\infty}$ are rotation
invariant, we have $$c\left(\sigma_{n,\,\lambda},H_{\varphi_{N}},H^{\infty}\right)=c\left(\sigma_{\mu,n},H_{\varphi_{N}},H^{\infty}\right)$$
for every $\lambda,\mu$ with $\vert\lambda\vert=\vert\mu\vert=r$.
Let $\lambda=-r$. To get a lower estimate for $\Vert\Psi\Vert_{H_{\varphi_{_{N}}}/b_{\lambda}^{n}H_{\varphi_{_{N}}}}$
consider $G$ such that $\Psi-G\in b_{\lambda}^{n}{\rm Hol}(\mathbb{D})$,
i.e. such that $bH_{n,\, N}\circ b_{\lambda}-G\circ b_{\lambda}\in z^{n}{\rm Hol}(\mathbb{D})$.

\vspace{0.2cm}
 3) First, we show that
\[
\psi=:\,\Psi\circ b_{\lambda}=\, bH_{n,\, N}\circ b_{\lambda}\]
is a polynomial (of degree $n$ if $\varphi=z$ and $2n$ if $\varphi=z^{2}$)
with positive coefficients. Note that
\[
Q_{n}\circ b_{\lambda}=\sum_{k=0}^{n-1}z^{k}\frac{(1-\vert\lambda\vert^{2})^{1/2}}{1-\overline{\lambda}b_{\lambda}(z)}=\left(1-\vert\lambda\vert^{2}\right)^{-\frac{1}{2}}\left(1+(1-\overline{\lambda})\sum_{k=1}^{n-1}z^{k}-\overline{\lambda}z^{n}\right)=\]
\[
=(1-r^{2})^{-1/2}\left(1+(1+r)\sum_{k=1}^{n-1}z^{k}+rz^{n}\right)=:(1-r^{2})^{-1/2}\psi_{1}.\]
Hence, $\psi=\Psi\circ b_{\lambda}=bH_{n,\, N}\circ b_{\lambda}=b\varphi_{N}\circ\left(\left(1-r^{2}\right)^{-\frac{1}{2}}\psi_{1}\right)$
and
\[
\varphi_{N}\circ\psi_{1}=\psi_{1}^{N}(z),\; N=1,\,2.\]

4) Next, we show that \[
\sum^{m}(\psi)=:\sum_{j=0}^{m}\hat{\psi}(j)\geq\left\{ \begin{array}{c}
(2\sqrt{2})^{-1}\sqrt{n(1-r)^{-1}}\;\;{\rm {\rm if}}\;\; N=1\\
16^{-1}n(1-r)^{-1}\;\;{\rm {\rm if}}\;\; N=2\end{array}\right.,\]
where $m=n/2$ if $n$ is even and $m=(n+1)/2$ if $n$ is odd.

\vspace{0.3cm}
Indeed, setting $S_{n}=\sum_{j=0}^{n}z^{j}$, we have both for $N=1$
and $N=2$\[
\sum^{m}\left(\psi_{1}^{N}\right)=\sum^{m}\left(\left(1+(1+r)\sum_{t=1}^{n-1}z^{t}+rz^{n}\right)^{N}\right)\geq\sum^{m}\left(S_{n-1}^{N}\right).\]
Next, we obtain\[
\sum^{m}\left(S_{n-1}^{N}\right)=\sum^{m}\left(\left(\frac{1-z^{n}}{1-z}\right)^{N}\right)=\]
\[
=\sum^{m}\left((1-z)^{-N}\right)=\frac{1}{(N-1)!}\sum^{m}\left(\frac{d^{N-1}}{dz^{N-1}}\frac{1}{1-z}\right)=\sum_{j=0}^{m}C_{N+j-1}^{j}.\]
Now if $N=1,$ then \[
\sum_{j=0}^{m}C_{N+j-1}^{j}=m+1\geq\frac{n}{2},\]
 whereas if $N=2$ then \[
\sum_{j=0}^{m}C_{N+j-1}^{j}=\frac{(m+1)(m+2)}{2}\geq\frac{(n+2)(n+4)}{8}\geq\frac{n^{2}}{8}.\]
Finally, since $\sum^{m}(\psi)=b\sum^{m}(\varphi_{N}\circ\psi_{1})=b\left(1-r^{2}\right){}^{-N/2}\sum^{m}(\psi_{1}^{N})$
we get\[
\sum^{m}(\psi)\geq\left\{ \begin{array}{c}
(2(1-r))^{-1/2}nb/2\;\;{\rm {\rm if}}\;\; N=1\\
(2(1-r))^{-1}n^{2}b/8\;\;{\rm {\rm if}}\;\; N=2\end{array}\right.,\]
with $b=\varphi_{N}(n)=\left\{ \begin{array}{c}
n^{-1/2}\;\;{\rm {\rm if}}\;\; N=1\\
n^{-1}\;\;{\rm {\rm if}}\;\; N=2\end{array}\right..$ This gives the result claimed.
 5) Let $F_{n}=\Phi_{m}+z^{m}\Phi_{m}$, where $\Phi_{k}$ stands
for the $k$-th Fejer kernel. We have $\Vert g\Vert_{\infty}\Vert F_{n}\Vert_{L^{1}}\geq\Vert g* F_{n}\Vert_{\infty}$
for every $g\in L^{\infty}\left(\mathbb{T}\right),$ and taking the
infimum over all $g\in H^{\infty}$ satisfying $\hat{g}(k)=\hat{\psi}(k),\;\forall k\in[0,\, n-1],$
we obtain\[
\Vert\psi\Vert_{H^{\infty}/z^{n}H^{\infty}}\geq\frac{{\displaystyle 1}}{{\displaystyle 2}}\Vert\psi*F_{n}\Vert_{\infty},\]
where $*$ stands for the usual convolution product. Now using part
4),
\[
\Vert\Psi\Vert_{H^{\infty}/b_{\lambda}^{n}H^{\infty}}=\Vert\psi\Vert_{H^{\infty}/z^{n}H^{\infty}}\geq\frac{{\displaystyle 1}}{{\displaystyle 2}}\Vert\psi* F_{n}\Vert_{\infty}\geq\]
\[
\geq\frac{{\displaystyle 1}}{{\displaystyle 2}}\left|\left(\psi*F_{n}\right)(1)\right|\geq\frac{{\displaystyle 1}}{{\displaystyle 2}}{\displaystyle \sum_{j=0}^{m}\hat{\psi}(j)}{\displaystyle \geq}\]
\[
\geq\left\{ \begin{array}{c}
(4\sqrt{2})^{-1}\sqrt{n(1-r)^{-1}}\;\;{\rm {\rm if}}\;\; N=1\\
32^{-1}n(1-r)^{-1}\;\;{\rm {\rm if}}\;\; N=2\end{array}\right..\]
6) In order to conclude, it remains to use (7.1.1). \end{proof}

\subsection{The case $X=H^{p},\,1\leq p\leq+\infty$}\label{subsection72}

Here we prove the sharpness (for even $p$) of the upper estimate
found in Theorem \ref{Theorem C}. We first prove the following lemma.

\begin{lem}\label{lemma721} Let $p$, $q$
such that $\frac{p}{q}\in\mathbb{Z}_{+}$, then $c\left(\sigma,\, H^{p},H^{\infty}\right)\geq c\left(\sigma,\, H^{q},H^{\infty}\right)^{\frac{q}{p}}$
for every sequence $\sigma$ of $\mathbb{D}$. \end{lem}
\begin{proof}
\textbf{Step 1.} Recalling that \[
c\left(\sigma,\, H^{p},\, H^{\infty}\right)=\sup_{\parallel f\parallel_{p}\leq1}\inf\left\{ \left\Vert g\right\Vert _{\infty}:\, g\in Y,\, g_{|\sigma}=f_{|\sigma}\right\} ,\]
we first prove that\[
c\left(\sigma,\, H^{p},\, H^{\infty}\right)=\sup_{\parallel f\parallel_{p}\leq1,\, f\, outer\,}\inf\left\{ \left\Vert g\right\Vert _{\infty}:\, g\in Y,\, g_{|\sigma}=f_{|\sigma}\right\} .\]
Indeed, we clearly have the inequality\[
\sup_{\parallel f\parallel_{p}\leq1,\, f\, outer\,}\inf\left\{ \left\Vert g\right\Vert _{\infty}:\, g\in Y,\, g_{|\sigma}=f_{|\sigma}\right\} \leq c\left(\sigma,\, H^{p},\, H^{\infty}\right),\]
and if the inequality were strict, that is to say \[
\sup_{\parallel f\parallel_{p}\leq1,\, f\, outer\,}\inf\left\{ \left\Vert g\right\Vert _{\infty}:\, g\in Y,\, g_{|\sigma}=f_{|\sigma}\right\} <\sup_{\parallel f\parallel_{p}\leq1}\inf\left\{ \left\Vert g\right\Vert _{\infty}:\, g\in Y,\, g_{|\sigma}=f_{|\sigma}\right\} ,\]
then we could write that there exists $\epsilon>0$ such that for
every $f=f_{i}.f_{o}\in H^{p}$ (where $f_{i}$ stands for the inner
function corresponding to $f$ and $f_{o}$ to the outer one) with
$\left\Vert f\right\Vert _{p}\leq1$ (which also implies that $\left\Vert f_{o}\right\Vert _{p}\leq1$,
since $\left\Vert f_{o}\right\Vert _{p}=\left\Vert f\right\Vert _{p})$,
there exists a function $g\in H^{\infty}$ verifying both $\left\Vert g\right\Vert _{\infty}\leq(1-\epsilon)c\left(\sigma,\, H^{p},\, H^{\infty}\right)$
and $g_{\vert\sigma}=f_{o\vert\sigma}.$ This entails that $f_{\vert\sigma}=\left(f_{i}g\right)_{\vert\sigma}$
and since $\left\Vert f_{i}g\right\Vert _{\infty}=\left\Vert g\right\Vert _{\infty}\leq(1-\epsilon)c\left(\sigma,\, H^{p},\, H^{\infty}\right)$
, we get that $c\left(\sigma,\, H^{p},\, H^{\infty}\right)\leq(1-\epsilon)c\left(\sigma,\, H^{p},\, H^{\infty}\right)$,
which is a contradiction and proves the equality of Step 1.

\textbf{Step 2.} Using the result of Step 1, we get that $\forall\epsilon>0$
there exists an outer function $f_{o}\in H^{q}$ with $\left\Vert f_{o}\right\Vert _{q}\leq1$
and such that\[
\inf\left\{ \left\Vert g\right\Vert _{\infty}:\, g\in Y,\, g_{|\sigma}=f_{o|\sigma}\right\} \geq c\left(\sigma,\, H^{q},H^{\infty}\right)-\epsilon.\]
Now let $F=f_{o}^{\frac{q}{p}}\in H^{p}$, then $\left\Vert F\right\Vert _{p}^{p}=\left\Vert f_{o}\right\Vert _{q}^{q}\leq1$.
We suppose that there exists $g\in H^{\infty}$ such that $g_{\vert\sigma}=F_{\vert\sigma}$
with \[
\left\Vert g\right\Vert _{\infty}<\left(c\left(\sigma,\, H^{q},H^{\infty}\right)-\epsilon\right)^{\frac{q}{p}}.\]
 Then, since $g\left(\lambda_{i}\right)=F\left(\lambda_{i}\right)=f_{o}\left(\lambda_{i}\right)^{\frac{q}{p}}$
for all $i=1,\,\dots,\, n$, we have $g\left(\lambda_{i}\right)^{\frac{p}{q}}=f_{o}\left(\lambda_{i}\right)$
and $g^{\frac{p}{q}}\in H^{\infty}$ since $\frac{p}{q}\in\mathbb{Z}_{+}$.
We also have \[
\left\Vert g^{\frac{p}{q}}\right\Vert _{\infty}=\left\Vert g\right\Vert _{\infty}^{\frac{p}{q}}<\left(c\left(\sigma,\, H^{q},H^{\infty}\right)-\epsilon\right)^{\frac{q}{p}},\]
which is a contradiction. As a result, we have\[
\left\Vert g\right\Vert _{\infty}\geq\left(c\left(\sigma,\, H^{q},H^{\infty}\right)-\epsilon\right)^{\frac{q}{p}},\]
for all $g\in H^{\infty}$ such that $g_{\vert\sigma}=F_{\vert\sigma}$,
which gives \[
c\left(\sigma,\, H^{p},H^{\infty}\right)\geq\left(c\left(\sigma,\, H^{q},H^{\infty}\right)-\epsilon\right)^{\frac{q}{p}},\]
and since that inequality is true for every $\epsilon>0$, we get
the result.
\end{proof}
\begin{proof}[Proof of Theorem \ref{Theorem C}  (the left-hand side inequality for $p\in2\mathbb{N}$, $p>2$   only)] We first prove the lower estimate for $c\left(\sigma_{n,\,\lambda},\, H^{p},H^{\infty}\right)\,.$
Writing $p=2(p/2)$, we apply Lemma \ref{lemma721} with $q=2$ and this gives\[
c\left(\sigma_{n,\,\lambda},\, H^{p},H^{\infty}\right)\geq c\left(\sigma_{n,\,\lambda},\, H^{2},H^{\infty}\right)^{\frac{2}{p}}\geq32^{-\frac{1}{p}}\left(n(1-\left|\lambda\right|)^{-1}\right)^{\frac{2}{p}}\]
for all integer $n\geq1$. The last inequality is a consequence of
Theorem \ref{Theorem C}  (left-hand side inequality) for the particular case $p=2$
which has been proved in Subsection 7.1. \end{proof}

\section{Comparing our results with Carleson interpolation}\label{section8}

Recall that given a (finite) set $\sigma=\{\lambda_{1},\,...,\,\lambda_{n}\}\subset\mathbb{D}$,
the Carleson interpolation constant $C_{I}(\sigma)$ is defined by\[
C_{I}(\sigma)=\sup_{\parallel a\parallel_{l^{\infty}}\leq1}\inf\left(\parallel g\parallel_{\infty}:\, g\in H^{\infty},\, g_{|\sigma}=a\right).\]
We introduce the evaluation functionals $\varphi_{\lambda}$ for $\lambda\in\mathbb{D}$,
as well as the evaluation of the derivatives $\varphi_{\lambda,s}$
($s=0,1,...)$ \[
\varphi_{\lambda}(f)=f(\lambda),\:\: f\in X,\,\mbox{\,\ and}\,\,\varphi_{\lambda,s}(f)=f^{(s)}(\lambda),\:\: f\in X.\]

\begin{thm} \label{Theorem E}Let $X$ be a Banach space,
$X\subset{\rm Hol}(\mathbb{D})$, and $\sigma=\{\lambda_{1},\,...,\,\lambda_{n}\}$
be a sequence of distinct points in the unit disc $\mathbb{D}.$ We
have,
\[
\max_{1\leq i\leq n}\left\Vert \varphi_{\lambda_{i}}\right\Vert \leq c(\sigma,\, X,\, H^{\infty})\leq C_{I}(\sigma)\max_{1\leq i\leq n}\left\Vert \varphi_{\lambda_{i}}\right\Vert ,\]
 where $C_{I}(\sigma)$ stands for the Carleson interpolation constant.
\end{thm}

Theorem \ref{Theorem E}  (see \cite{Z3} for its proof) tells us that, for $\sigma$
with a {}``reasonable'' interpolation constant $C_{I}(\sigma)$,
the quantity $c(\sigma,\, X,\, H^{\infty})$ behaves as $\max_{i}\left\Vert \varphi_{\lambda_{i}}\right\Vert $.
However, for {}``tight'' sequences $\sigma$, the constant $C_{I}(\sigma)$
is so large that the estimate in question contains almost no information.
On the other hand, an advantage of the estimate of Theorem \ref{Theorem E}  is that
it does not contain ${\rm card\,}\sigma=n$ explicitly. Therefore,
for well-separated sequences $\sigma,$ Theorem \ref{Theorem E}  should give a better
estimate than those of Theorem \ref{Theorem C}  and Theorem \ref{Theorem D} .

Now, how does the interpolation constant $C_{I}(\sigma)$ behave in
terms of the characteristics $r$ and $n$ of $\sigma$? We answer
this question in \cite{Z3} for some particular sequences $\sigma$.
More precisely, we compare these quantities for the cases $X=H^{2}$,
$X=L_{a}^{2}$ and for three geometrically simple configurations:
\textit{two-points sets} $\sigma$, \textit{circular} and \textit{radial}
sequences $\sigma$.

Let us recall that our specific upper bounds in Theorems C and D are
sharp over all $n$ elements sequences $\sigma.$ However, we give
in \cite{Z3} some very special \textit{radial} and \textit{circular}
sequences $\sigma$ such that the estimate of $c(\sigma,\, H^{2},\, H^{\infty})$
via the Carleson constant $C_{I}(\sigma)$ (using Theorem \ref{Theorem E} ) is comparable
with or better than the estimates of Theorem \ref{Theorem C}  (for $X=H^{2})$ and
Theorem \ref{Theorem D}  (for $X=L_{a}^{2}$) . We also give some examples of \textit{radial}
and \textit{circular} sequences but also of \textit{two-points sets},
such that it is worse (i.e. for which our estimate is better). More
specific \textit{radial sequences} are studied in \cite{Z3}:\textit{
sparse sequences, condensed sequences }and\textit{ long sequences. }

{\small \hfill Received July 4, 2010}

{\small \hfill Revised version received October 25, 2010}

\end{document}